\documentclass{jgmOF} 
\usepackage{amsmath}
\usepackage{paralist}
\usepackage{epsfig} 
\usepackage{epstopdf} 

\allowdisplaybreaks

\def\currentvolume{X}
 \def\currentissue{X}
  \def\currentyear{200X}
   \def\currentmonth{XX}
    \def\ppages{X--XX}
     \def\DOI{10.3934/jgm.2022014}

\usepackage[noabbrev,capitalise]{cleveref}
\usepackage{autobreak}
\usepackage{subcaption}
\usepackage{forest}
\forestset{
	default preamble={
		for tree={fill,circle,draw,grow=east,for descendants={fill=white}}
	}
}
\usepackage{scalerel}
\usepackage{tikz}
\usetikzlibrary{svg.path}

\definecolor{orcidlogocol}{HTML}{A6CE39}
\tikzset{
  orcidlogo/.pic={
    \fill[orcidlogocol] svg{M256,128c0,70.7-57.3,128-128,128C57.3,256,0,198.7,0,128C0,57.3,57.3,0,128,0C198.7,0,256,57.3,256,128z};
    \fill[white] svg{M86.3,186.2H70.9V79.1h15.4v48.4V186.2z}
                 svg{M108.9,79.1h41.6c39.6,0,57,28.3,57,53.6c0,27.5-21.5,53.6-56.8,53.6h-41.8V79.1z M124.3,172.4h24.5c34.9,0,42.9-26.5,42.9-39.7c0-21.5-13.7-39.7-43.7-39.7h-23.7V172.4z}
                 svg{M88.7,56.8c0,5.5-4.5,10.1-10.1,10.1c-5.6,0-10.1-4.6-10.1-10.1c0-5.6,4.5-10.1,10.1-10.1C84.2,46.7,88.7,51.3,88.7,56.8z};
  }
}

\newcommand\orcidicon[1]{\href{https://orcid.org/#1}{\mbox{\scalerel*{
\begin{tikzpicture}[yscale=-1,transform shape]
\pic{orcidlogo};
\end{tikzpicture}
}{|}}}}

\def\rg{\mathrm{rg }\,  }
\def\d{\mathrm{d}}
\def\p{\partial }

\def\D{\mathrm{D}}

\def\Z{\mathbb{Z}}

\def\N{\mathbb{N}}
\def\R{\mathbb{R}}

\newtheorem{assumption}{Assumption}[section]

\newtheorem{theorem}{Theorem}[section]
\newtheorem{corollary}{Corollary}

\newtheorem{lemma}[theorem]{Lemma}
\newtheorem{proposition}{Proposition}

\theoremstyle{definition}

\newtheorem{remark}{Remark}

\title[Backward Error Analysis for Variational Discretisations of PDEs]
{Backward Error Analysis for Variational Discretisations of PDEs} 

\author[Robert I McLachlan and Christian Offen]{}


 \email{R.McLachlan@massey.ac.nz}
 \email{christian.offen@uni-paderborn.de}

\thanks{$^*$Corresponding author: Christian Offen $^{\protect \orcidicon{0000-0002-5940-8057}}$}

\begin{document}
\maketitle

\centerline{\scshape Robert I McLachlan}
\medskip
{\footnotesize
 \centerline{Massey University}
   \centerline{Private Bag 11 222, Palmerston North, 4442}
   \centerline{New Zealand}
} 

\medskip

\centerline{\scshape Christian Offen$^*$ $^{\protect \orcidicon{0000-0002-5940-8057}}$}
\medskip
{\footnotesize
 \centerline{Paderborn University}
   \centerline{Warburger Str. 100, 33098 Paderborn}
   \centerline{Germany}
}

\bigskip

 \centerline{(Communicated by the associate editor name)}


\begin{abstract}
In backward error analysis, an approximate solution to an equation is compared to the exact solution to a nearby `modified' equation. In numerical ordinary differential equations, the two agree up to any power of the step size. If the differential equation has a geometric property then the modified equation may share it. In this way, known properties of differential equations can be applied to the approximation. But for partial differential equations, the known modified equations are of higher order,  limiting applicability of the theory. Therefore, we study symmetric solutions of discretized partial differential equations that arise from a discrete variational
principle. These symmetric solutions obey infinite-dimensional functional
equations. We show that these equations admit second-order modified
equations which are Hamiltonian and also possess first-order Lagrangians in
modified coordinates. The modified equation and its associated structures
are computed explicitly for the case of rotating travelling waves
in the nonlinear wave equation.
\end{abstract}



\section{Introduction}
\label{sec:introduction}
\subsection*{Motivation}
Backward error analysis is a key tool for understanding the behaviour of discretizations of differential equations. The numerical solution of an ODE or PDE closely approximates the exact solution of a {\em modified equation} at the grid points. This approximation is not exact, but can typically be made of
arbitrarily high order in the step sizes.
The modified equation  can be obtained as a series expansion of the discretisation in the step sizes. Finding the modified equation and analysing its properties is referred to as {\em backward error analysis (BEA)} (see, for instance \cite[\S IX]{GeomIntegration}).

For Hamiltonian ODEs discretised by a symplectic integrator the modified equation is itself Hamiltonian with respect to a modified Hamiltonian and the original symplectic structure and phase space. The Hamiltonian is given as a formal power series in the time step which typically does not converge. However, optimal truncation results are available \cite{GeomIntegration}. A variational version of backward error analysis was developed in \cite{Vermeeren2017}. Next to the analysis of numerical methods for ODEs, backward error analysis has been employed to improve machine learning techniques \cite{LagrangianShadowIntegrators,symplecticShadowIntegrators}.

Backward error analysis for Hamiltonian PDEs has been developed by Moore and Reich \cite{moore2003backward} and by Islas and Schober \cite{islas2005backward}.
We illustrate this briefly on the main example of the paper, the nonlinear wave equation
\begin{equation}
	\label{eq:NonLinearWavePDE}
	u_{tt}-u_{xx}-\nabla W(u)=0,\quad u\colon\R^2\to\R^d,\ W\colon\R^d\to\R
\end{equation}
and its five-point stencil discretisation
\begin{equation}\label{eq:5ptStencilOnNonWaveEQ}
	\begin{aligned}
		0&=\frac{1}{{\Delta t}^2} \left( u(t-\Delta t, x) - 2 u(t,x) + u(t+\Delta t,x) \right) \\
		&- \frac{1}{{\Delta x}^2} \left( u(t, x - \Delta x) - 2 u(t,x) + u(t,x+\Delta x) \right) \\
		&-\nabla W(u(t,x)).
	\end{aligned}
\end{equation}
Expanding \eqref{eq:5ptStencilOnNonWaveEQ} as a formal power series in $\Delta t$ and $\Delta x$ gives a modified equation
\begin{equation}
	\label{eq:mod}
	0 = u_{tt} - u_{xx} - \nabla W(u) + \frac{\Delta t^2}{12} u_{tttt} - \frac{\Delta x^2}{12}u_{xxxx},
\end{equation}
where terms of order 4 in $\Delta t$ or $\Delta x$ have been truncated.
This modified equation does preserve some features of the original
equation (\ref{eq:NonLinearWavePDE}). It is variational and multisymplectic. It has variational
symmetries (translations in $x$ and in $t$) that can be used to construct approximate conservation
laws of (\ref{eq:5ptStencilOnNonWaveEQ}). However, it is of higher
order than (\ref{eq:NonLinearWavePDE}). Its multisymplectic formulation needs more dependent variables, while its Lagrangian formulation is second rather than first order.
In contrast to the ODE case, the higher derivatives
cannot be  eliminated. Thus, for PDEs a key advantage of backward error analysis for ODEs---that the modified equation lies in the same class as the original---is lost. Moreover, optimal truncation techniques have not yet been developed \cite[\S 5.3.2]{MarsdenWestVariationalIntegrators}.

In this paper we study backward error analysis for methods such as (\ref{eq:5ptStencilOnNonWaveEQ}) through the lens of {\em symmetric solutions}.
We restrict our attention to an analysis of the structure of the modified equation as a formal power series in the discretisation parameters without discussing convergence issues.
Clearly, solutions of (\ref{eq:5ptStencilOnNonWaveEQ}) that are independent of $x$ or $t$ have a standard Hamiltonian modified ODE. Travelling wave solutions of the form $u(x,t)=\phi(\xi)$, $\xi := x-ct$ are invariant under the symmetry with generator
$\partial_t + c\partial_x$. They obey the discrete travelling wave equation
\begin{equation}
	\label{eq:dtwe}
	\begin{aligned}
		0&=\frac{1}{{\Delta t}^2} \left( \phi(\xi+c\Delta t) - 2 \phi(\xi) + \phi(\xi-c\Delta t) \right) \\
		&- \frac{1}{{\Delta x}^2} \left(  \phi(\xi+\Delta x) - 2 \phi(\xi) + \phi(\xi-\Delta x) \right) \\
		&-\nabla W(\phi(\xi)).
	\end{aligned}
\end{equation}

It was shown in \cite{mcdonald2016,mcdonald2013travelling} that this essentially infinite-dimensional nonlinear functional equation has a second-order Hamiltonian modified equation. Its Hamiltonian and symplectic structure were computed for $\phi\in\mathbb{R}$ up to fourth order. This
reduction in phase space dimension (from infinity to two) is a particularly
dramatic example of backward error analysis, and motivates us to extend this example to a wide
class of discrete methods and symmetries.

\begin{remark}
	Equation \eqref{eq:5ptStencilOnNonWaveEQ} has continuous independent variables $(t,x)\in\mathbb{R}^2$.
	Solutions restricted to the grid $\Delta t \mathbb Z\times\Delta x\mathbb{Z}$ satisfy the standard
	five point stencil. Here $\Z$ denotes the set of integers.
	We adopt this point of view because it removes the awkward distinction
	between the discrete translation symmetry of the grid and the continuous translation symmetry
	of the PDE, and because \eqref{eq:dtwe} necessarily has a continuous independent variable $\xi\in\mathbb{R}$. Moreover, we restrict our attention to an analysis of the formal structure of \eqref{eq:5ptStencilOnNonWaveEQ} rather than developing a functional analytic setting and considering boundary conditions.
\end{remark}

\begin{remark}
	The PDE \eqref{eq:NonLinearWavePDE} has another symmetry, the hyperbolic rotation
	with generator $t\partial_x + x\partial_t$, which is not shared by the discretisation \eqref{eq:5ptStencilOnNonWaveEQ}. In the following we will restrict our discussions to symmetries that are exactly preserved by the discretisation. We will consider in detail the case in which $W$ is of the form $W = \frac 12 V(\|u\|^2)$ and $u$ is $\R^2$-valued, leading to the symmetry group with generators $\partial_x$, $\partial_t$, and $u_2\partial_{u_1} - u_1\partial_{u_2}$ and to rotating travelling waves.
\end{remark}

\subsection*{Palais' principle of symmetric criticality}
There is an extensive and well-known theory of group-invariant solutions of partial differential equations
\cite{olver1986}. We are particularly interested in cases that reduce to an ordinary differential equation.
Many partial differential equations fulfil a variational principle, i.e.\ they arise as the Euler--Lagrange equations corresponding to an action functional $S\colon U \to \R$ of the form
\begin{equation}\label{eq:IntroS}
	S(u) = \int \mathbb L (t_1,\ldots,t_n,u,u_{t_j},u_{t_i t_j},\ldots)\, \d t_1 \ldots \d t_n
\end{equation}
defined on a suitable function space $U$ (typically a Banach space). Here, the independent variables are denoted by $\mathbf{t}=(t_1,\ldots,t_n)$ and can refer to space and time dimensions.
However,  the reduced equations of the group-invariant solutions of a variational PDE are not
necessarily  variational. (This occurs even for standard examples in general relativity
\cite{Fels2002}.) This situation is the subject of {\em Palais' principle of symmetric criticality} \cite{palais1979}, which is formulated for general functionals $S\colon U \to \R$ (not necessarily of the form \eqref{eq:IntroS}).

Consider the action of a Lie group $G$ on a function space $U$. Let us denote the set of elements $u \in U$ which are fixed points under the action of the symmetry group $G$ by $U^{\mathrm{sym}}$, i.e.\ $U^{\mathrm{sym}} = \{ u \in U \, | \, g \cdot u = u \; \forall g \in G\}$.
Assume that $U^{\mathrm{sym}}$ is a submanifold of $U$.
Critical points of the action $S\colon U \to \R$ which lie in $U^{\mathrm{sym}}$ are critical points of the restricted functional $S|_{U^{\mathrm{sym}}}\colon U^{\mathrm{sym}} \to \R$. If the converse holds true as well, i.e.\ if the critical points of $S|_{U^{\mathrm{sym}}}$ are critical points of $S$, then we say the {\em principle of symmetric criticality} holds true. In other words, the principle of symmetric criticality says that symmetric elements $u \in U^{\mathrm{sym}}$ which are stationary points of $S$ with respect to symmetric variations are stationary with respect to all variations.
Palais analyses in \cite{palais1979} when the principle of symmetric criticality applies. He proves in particular that the principle holds if the symmetry group is compact or the group action is isometric and $U$ is a Banach space.

We will restrict attention to cases where the principle of symmetric criticality applies. This is
easy to check in specific examples.

\subsection*{Variational structure of symmetric solutions of discrete systems}

Variational principles are also useful for constructing numerical integration schemes.
In the ODE case, discrete variational integrators are automatically symplectic and show excellent energy conservation properties as well as favourable preservation properties of the topological structures of phase portraits when applied to Hamiltonian systems. Moreover, the discretised variational principle allows for a theoretical analysis using discrete versions of tools known from the continuous setting such as, for example, the discrete Noether theorem. In the PDE case, discrete
variational integrators obey a discrete multisymplectic conservation law  \cite{MarsdenWestVariationalIntegrators}.

The method \eqref{eq:5ptStencilOnNonWaveEQ} has Lagrangian
\[
L_\Delta = \frac{\|u(t-\Delta t,x)-u(t,x)\|^2}{2\Delta t^2} - \frac{\|u(t,x-\Delta x)-u(t,x)\|^2}{2\Delta x^2}-W(u(t,x))
\]
which, restricted to the grid, becomes the standard discrete Lagrangian of the five point stencil. It approximates the Lagrangian $\mathbb{L}(u,u_t,u_x) = \frac 12 (\|u_t\|^2-\|u_x\|^2)-W(u)$.

\begin{figure}
	\begin{center}
		\includegraphics[width=0.8\textwidth]{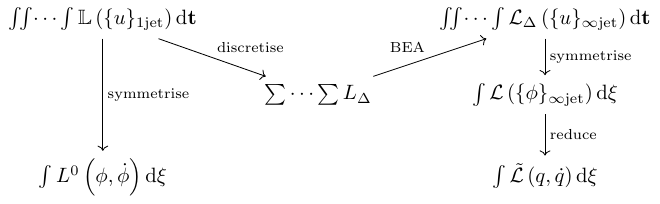}
	\end{center}
	\caption[Setting theorem]{Illustration of \cref{thm:ConsPalais}.
		The left hand column gives the actions of a PDE and an associated ODE that governs its symmetric solutions such as travelling waves.
		The right hand column gives three Lagrangians of modified equations of a variational discretization.
		Top: of the discretization, containing arbitrarily high derivatives;
		middle: of its symmetric solutions, containing arbitrarily high derivatives; and
		bottom: of its symmetric solutions, containing first derivatives only. $\tilde{\mathcal L}$ can be regarded
		as a modified Lagrangian of $L^0$.
		${\mathcal L}_\Delta$, $\mathcal L$ and $\tilde {\mathcal L}$ are formal power series in the step sizes.}\label{fig:SettingThm}
\end{figure}

Our main result, stated in Theorem \ref{thm:ConsPalais} and illustrated in Figure \ref{fig:SettingThm}, uses a blend of Hamiltonian and Lagrangian methods. It uses similar ideas as our discussion of symmetric linear multistep methods in \cite{BEAMulti}.
Essentially, it is easy to check (at least in examples) that discrete symmetric solutions (such as those obeying \eqref{eq:dtwe}) obey a second order modified equation. The theorem states that this modified equation is also variational, and
satisfies a standard first order variational principle in a sense to be made precise.

If $u$ is a function of independent variables $\mathbf{t}=(t_1,\ldots,t_n)$, we denote the jet of order $K$ of $u$ by $\{u\}_{K\mathrm{jet}}$, where $K \in \N \cup \{\infty\}$.

\begin{theorem}\label{thm:ConsPalais}
	Consider a first-order Lagrangian density $\mathbb{L}(\{u\}_{1\mathrm{jet}})\d \mathbf{t}$, where $\mathbf{t}=(t_1,\ldots,t_n)$ denotes the independent variables, and a consistent discrete Lagrangian $L_\Delta$.
	A series expansion of the discrete Lagrangian $L_\Delta$ in the step sizes yields a Lagrangian density $\mathcal{L}_\Delta(\{u\}_{\infty\mathrm{jet}})\d \mathbf{t}$ of infinite order given as a formal power series. Let $\mathcal E$ be the Euler--Lagrange operator. Consider a symmetry action such that $\mathcal{L}_\Delta(\{u\}_{\infty\mathrm{jet}})\d \mathbf{t}$ reduces to a Lagrangian density of the form $\mathcal{L}(\{\phi\}_{\infty\mathrm{jet}})\d \xi$, where $\xi$ is 1-dimensional. Scaling all step sizes by a formal variable $h$, assume that $\mathcal{E} \mathcal{L}=0$ is equivalent to the power series
	\begin{equation}\label{eq:ODEThm}
		\ddot \phi = a_0(\phi,\dot \phi) + \sum_{i=1}^\infty h^i a_i(\{\phi\}_{m_i\mathrm{jet}}).
	\end{equation}
	
	If \eqref{eq:ODEThm} is equivalent to reducing the Euler--Lagrange equations $\mathcal E \mathbb{L}(\{u\}_{1\mathrm{jet}}) = 0$ by the considered symmetry (i.e.\ Palais' principle of symmetric criticality holds), then, under non-degeneracy assumptions, \eqref{eq:ODEThm} is formally governed by the Euler--Lagrange equations of a first-order Lagrangian $\tilde{\mathcal L}(q,\dot q)\, \d \xi$.
\end{theorem}

In other words, the theorem says that modified variational principles corresponding to variational integrators for symmetric solutions have variational structure of the correct order.

\begin{remark}
	The condition that $\mathcal{E} \mathcal{L}=0$ is equivalent to the power series \eqref{eq:ODEThm}, is equivalent to the condition that the $h^0$-term $L^0$ of the power series of the Lagrangian $\mathcal L$ constitutes a non-degenerate Lagrangian, i.e., that the matrix $\left(\frac{\p^2 L^0}{\p \phi_i \p \phi_j}\right)_{i,j}$ is invertible. Indeed, as the discretisation is consistent and Palais' principle of symmetric criticality is assumed to hold, this term coincides with  the symmetrisation $L^0$ of $\mathbb L$ in \cref{fig:SettingThm}. The condition then says that the Lagrangian $L^0$ is non-degenerate.
\end{remark}

\subsection*{Choice of coordinates}
Using the notation of \cref{fig:SettingThm}, notice that the modified Lagrangian $\tilde{\mathcal L}$ needs to be expressed in new coordinates $(q,\dot q)$. We will see that if it is expressed in $\phi$ and its derivatives then it might contain second derivatives of $\phi$. Its Euler--Lagrange equations yield a jet-extension (i.e., a prolongation) of \eqref{eq:ODEThm}.
A change of coordinates from the variables $(\phi,\dot \phi)$ of the continuous Lagrangian $L$ to the variables $(q,\dot q)$ of the modified Lagrangian $\tilde{\mathcal L}$ does not in general  admit a closed form. However, we will show that there does exist an explicit description of the first-order system as a Hamiltonian system with a modified symplectic structure. Moreover, we will give a sufficient criterion for
the existence of a modified first-order Lagrangian in the original variables $(\phi,\dot \phi)$. Notice that in the literature the statement {\em has variational structure} sometimes requires the Lagrangian to be of the required order in the original variables, see \cite{Barbero2018}, for instance.

Additionally, we will verify that constants of motion are preserved when reducing the order of the Lagrangian.

\subsection*{Structure of the paper}
The remainder of the paper is structured as follows. In \cref{sec:GerneralTreatment} we present the method of reducing high-order Lagrangians that have a series structure. This provides a proof of \cref{thm:ConsPalais}. Moreover, we prove that conserved quantities and symmetries are passed on to the reduced system.

In \cref{sec:TravWave} we carry out the constructions of the proof of \cref{thm:ConsPalais} in detail for the example of the nonlinear wave equation \eqref{eq:NonLinearWavePDE} and its five-point stencil  \eqref{eq:5ptStencilOnNonWaveEQ}. We also consider the case in which $W$ is of the form $W = \frac 12 V(\|u\|^2)$ and $u$ is $\R^2$-valued, for which the continuous and discrete equations both admit rotating travelling waves $u(t,x)=R(t)\phi(x-ct)$, where $R$ is a rotation matrix. Such
waves satisfy a 4-dimensional Hamiltonian system. We
compute the modified Hamiltonian and symplectic structure of the discrete rotating travelling waves. The modified Lagrangian is computed for special cases.

A treatment of the special case of non-rotating travelling waves is contained in \cref{sec:MSApproach}, which shows an alternative theoretical approach using linear multistep methods and illustrates a relation to $P$-series. ($P$-series occur, for instance, in the analysis of partitioned Runge--Kutta methods. See \cite[III.2.1]{GeomIntegration} for an introduction.) The extra information is then used in a computational example to compute the modified Lagrangian efficiently. In addition, the modified Lagrangian is then of first order in the same variable $\phi$ as in the continuous setting. Similar results hold for rotating travelling waves in the case of a standing wave $c=0$ and when the step sizes fulfil the relation $\Delta x =c \Delta t$.


Source code for the computational examples of this work and a documentation of computational results can be found in \cite{multisymplecticSoftware}.

\section{Reduction of power series of high-order Lagrangians}\label{sec:GerneralTreatment}

Consider the following formal variational principle
\begin{equation}\label{eq:VarPrinc}
	\begin{split}
		\mathcal{S}(\phi)&=\int \mathcal{L}(\{\phi\}_{\infty\mathrm{jet}}(\xi))\d \xi \\
		&= \int \big( \mathcal{L}^0(\phi(\xi),\dot \phi(\xi)) + h \mathcal{L}^1(\{\phi\}_{M_1\mathrm{jet}}(\xi))+ h^2 \mathcal{L}^2(\{\phi\}_{M_2\mathrm{jet}}(\xi))+\ldots \big) \,\d \xi\\
		&=\int \sum_{i=0}^\infty h^i \mathcal{L}^i(\{\phi\}_{M_i\mathrm{jet}}(\xi)) \, \d \xi.
	\end{split}
\end{equation}
The Lagrangian $\mathcal{L}$ is given as a formal power series in the series parameter $h$. The expression $\{\phi\}_{M_i\mathrm{jet}}=\left(\phi,\dot \phi,\phi^{(2)},\ldots,\phi^{(M_i)}\right)$ denotes the jet of $\phi$ of order $M_i$. In the above expression, $\phi$ is a map that depends on a 1-dimensional variable $\xi$ and takes values in an $n$-dimensional manifold $Q$ which is locally identified with open subsets of $\R^n$. Our analysis focuses on local properties within a coordinate patch and does not consider global topological aspects.
In the following, we will also use $\phi$ to denote a variable on $Q$ and $\{\phi\}_{K\mathrm{jet}}$ to denote variables on the $K$-jet-space $\mathrm{Jet}_K(Q)$ over $Q$.
All maps are required to be sufficiently regular such that all considered derivatives exist.

We define the total derivative operator $\frac{\d}{\d \xi}\colon \mathcal C^\infty(\mathrm{Jet}_{K}(Q))\to  \mathcal C^\infty(\mathrm{Jet}_{K+1}(Q))$ acting on a function $\rho\colon \mathrm{Jet}_{K}(Q)\to \R$ in the variables of the jet space $\{\phi\}_{K\mathrm{jet}}=(\phi,\dot \phi,\phi^{(2)},\ldots,\phi^{(K)})$ as
\[
\frac{\d}{\d \xi} \rho\left(\{\phi\}_{K\mathrm{jet}}\right) = \sum_{i=0}^{K} \left\langle \nabla_{\phi^{(i)}}\rho, \phi^{(i+1)}\right\rangle.
\]
In other words, $\frac{\d}{\d \xi}$ acts on the expression $ \rho\left(\{\phi\}_{K\mathrm{jet}}\right)$ as if the $\phi$s were $\xi$-dependent functions.
The Euler--Lagrange equations to \eqref{eq:VarPrinc} are obtained by applying the Euler operator
\begin{equation}
	\begin{split}
		\mathcal{E}^j_\infty &= \frac{\p}{\p \phi_j}
		- \frac{\d}{\d \xi} \frac{\p}{\p \dot\phi_j}
		+\frac{\d^2}{\d \xi^2} \frac{\p}{\p \phi^{(2)}_j}+\ldots\\
		&= \sum_{i=0}^\infty (-1)^i \frac{\d^i}{\d \xi^i} \frac{\p}{\p \phi^{(i)}_j}
	\end{split}
\end{equation}
for each component $j=1,\ldots,n$ to the Lagrangian $\mathcal{L}$.
In the following $\mathcal{E}_\infty \mathcal{L}$ denotes the $n$-tuple $(\mathcal{E}^1_\infty \mathcal{L},\ldots, \mathcal{E}^n_\infty \mathcal{L})$. We also define
$\mathcal{E}^j_K =  \sum_{i=0}^K (-1)^i \frac{\d^i}{\d \xi^i} \frac{\p}{\p \phi^{(i)}_j}$ and $\mathcal{E}_K \mathcal{L} = (\mathcal{E}^1_K \mathcal{L},\ldots, \mathcal{E}^n_K \mathcal{L})$.
Under the non-degeneracy assumption that $\left(\frac{\p ^2\mathcal{L}^0}{\p \dot \phi_i \p \dot \phi_j}\right)_{0\le i,j \le n}$ is invertible, the Euler-Lagrange equations $\mathcal E_\infty(\mathcal L)=0$ yield the following ordinary differential equation given by a  formal power series.
\begin{equation}\label{eq:formalODE}
	\begin{split}
		\ddot \phi &= a_0(\phi,\dot \phi) + h a_1(\{\phi\}_{M_1\mathrm{jet}})+ h^2 a_2(\{\phi\}_{M_2\mathrm{jet}})+\ldots\\
		&= a_0(\phi,\dot \phi) + \sum_{i=1}^\infty h^i a_i(\{\phi\}_{M_i\mathrm{jet}})
	\end{split}
\end{equation}
The symbols $a_i$ denote sufficiently regular $\R^n$-valued maps.
We now claim that for any $N\in \N$ the truncation of \eqref{eq:formalODE} to order $N$, i.e.\
\begin{equation}\label{eq:truncatedODE}
	\ddot \phi = a_0(\phi,\dot \phi) + \sum_{i=1}^N h^i a_i(\{\phi\}_{M_i\mathrm{jet}})
\end{equation}
is formally $\mathcal{O}(h^{N+1})$-close to a second order equation
\begin{equation}\label{eq:2ndOrderODE}
	\ddot \phi = a_0(\phi,\dot \phi) + \sum_{i=1}^N h^i \widetilde{a_i}(\phi,\dot \phi).
\end{equation}
This can be seen by repeatedly substituting derivatives of order $j\ge2$ on the right hand side of \eqref{eq:truncatedODE} with the expression $\phi^{(j)}$ obtained by taking the $(j-2)$-derivative of \eqref{eq:truncatedODE} and truncating $\mathcal{O}(h^{N+1})$ terms.

Conserved quantities and symmetries of the original equation \eqref{eq:formalODE} are inherited by the reduced system \eqref{eq:2ndOrderODE} up to any order. This is made precise in the following proposition.

\begin{proposition}[Preservation of conserved quantities under reduction]\label{prop:presConsQ}
	Let $N \in \N$ be a truncation index and $M = \max\{M_1,\ldots,M_N\}$. If $I \colon \mathrm{Jet}_M(Q) \to \R$ is conserved along solutions of \eqref{eq:truncatedODE} up to $\mathcal O(h^{N+1})$ terms then $I$ induces a quantity $\overline I \colon \mathrm{Jet}_1(Q) \to \R$ that is conserved along solutions of the reduced system \eqref{eq:2ndOrderODE} up to $\mathcal O(h^{N+1})$ terms. If $\gamma$ is a solution of the truncated, reduced system \eqref{eq:2ndOrderODE} and $\mathrm{jet}_M(\gamma)$ is its prolongation then $\overline I \circ \gamma = I \circ \mathrm{jet}_M(\gamma) + \mathcal{O}(h^{N+1})$.
\end{proposition}

\begin{proof}
	To obtain $\overline I$ from $I$ we replace all 2nd and higher derivatives in $I$ by expressions in $(\phi,\dot \phi)$. These are obtained from (derivatives of) \eqref{eq:2ndOrderODE}. Substitutions may need to be iterated and order $\mathcal{O}(h^{N+1})$ terms are truncated. This corresponds to an {\em on-shell} expression of $I$.
	Prolongations of solutions to the truncated, reduced system \eqref{eq:2ndOrderODE} solve the truncated original system \eqref{eq:truncatedODE} up to higher order. Thus, a prolongation of a solution $\gamma$ to the truncated, reduced system conserves $I$ up to higher order. By construction of $\overline I$, it follows that $\overline I \circ \gamma = I \circ \mathrm{jet}_M(\gamma) + $ terms of higher order in $h$.
\end{proof}

\begin{remark}
	If $I$ is a conserved quantity in the sense of \cref{prop:presConsQ} of the high-order system and a conserved quantity $\overline I$ of the reduced system is constructed as above then $\overline I$ is a conserved quantity of the high-order system as well.
	The complement of the set described by \eqref{eq:truncatedODE} is open and dense in the jet-space of order $\max(M_1,\ldots,M_N)$ (unless degenerate). By only partially substituting higher order derivatives with expressions in lower order derivatives, one can obtain many conserved quantities of the high-order system, which will, unless degenerate, be functionally-independent on a dense open subset of the phase space.
\end{remark}

\begin{proposition}
	Consider a diffeomorphism $\chi \colon Q \to Q$ which is a symmetry of the action functional $\mathcal S$ from \eqref{eq:VarPrinc}, i.e.\
	\[\chi^\ast(  \mathcal{L}(\{\phi\}_{\infty\mathrm{jet}}(\xi))\d \xi) = \mathcal{L}(\{\phi\}_{\infty\mathrm{jet}}(\xi))\d \xi,\]
	where $\chi^\ast$ denotes the pullback map. Then $\chi$ is a symmetry of the reduced equation \eqref{eq:2ndOrderODE}, i.e.\ for any truncation index $N$, if $\phi$ solves \eqref{eq:2ndOrderODE} up to $\mathcal O(h^{N+1})$-terms then $\chi \circ \phi$ solves \eqref{eq:2ndOrderODE} up to $\mathcal O(h^{N+1})$-terms.
\end{proposition}

\begin{proof}
	If $\phi$ solves \eqref{eq:2ndOrderODE} then $\phi$ solves \eqref{eq:truncatedODE} up to higher order terms. As $\chi$ is a symmetry of $\mathcal S$, the curve $\chi \circ \phi$ solves \eqref{eq:truncatedODE} up to higher order terms. Thus, $\chi \circ \phi$ solves \eqref{eq:2ndOrderODE} up to higher order terms.
\end{proof}

Using the Legendre transformation $\mathcal H^0$ of $\mathcal L^0$ the Euler--Lagrange equations to zeroth order in $h$, i.e.\ $\mathcal E_1 \mathcal L^0 = 0$, can be transformed into a Hamiltonian system $(\mathrm{Jet}_1(Q),\omega^0,\mathcal H^0)$ with
\[
\omega^0 = \sum_{i=1}^n \d \mathfrak{q}^i \wedge \d \mathfrak{p}_i,
\quad \mathfrak{q} = \phi, \quad \mathfrak{p} = \frac{\p \mathcal L^0}{\p \dot \phi}
\]
and
\[
\mathcal H^0 = \langle \dot {\mathfrak{q}}, \mathfrak{p}\rangle_{\R^n} - \mathcal L^0.
\]
As we are assuming that the Lagrangian $\mathcal{L}^0$ is non-degenerate, i.e.\ $\left(\frac{\p ^2\mathcal{L}^0}{\p \dot \phi_i \p \dot \phi_j}\right)_{0\le i,j \le n}$ is of full rank, the Hamiltonian $\mathcal H^0$ can be expressed in the coordinates $\mathfrak{q},\mathfrak{p}$.

To formulate the main theorem of this section, we require that the truncation ${\mathcal L}^{[N]}$ of $\mathcal L$ after terms of order $N$ in $h$ constitutes a high order regular Lagrangian. This guarantees that a high order version of the Legendre transformation exists.

\begin{assumption}\label{ass:SNonDeg}
	For $N \in \N$ let $\mathcal{L}^{[N]}$ denote the truncation of the series $\mathcal{L}$ from \eqref{eq:VarPrinc} after $\mathcal{O}(h^N)$ terms. Define $M = \max\{M_1,\ldots,M_N\}$ as the order of the highest derivative that occurs in the expression $\mathcal{L}^{[N]}$.
	We assume that $\mathcal{L}^{[N]}$ constitutes a regular order $M$ Lagrangian, i.e.
	\[
	\mathrm{Hess}_{\phi^{(M)}}(\mathcal{L}^{[N]})= \left(\frac{\p ^2\mathcal{L}^{[N]}}{\p \phi^{(M)}_i \p \phi^{(M)}_j}\right)_{0\le i,j \le n}
	\]
	is of full rank.
	Similarly, we require the zeroth order Lagrangian to be non-degenerate, i.e.\ $\mathrm{Hess}_{\dot{\phi}}(\mathcal{L}^0)= \left(\frac{\p ^2\mathcal{L}^0}{\p \dot \phi_i \p \dot \phi_j}\right)_{0\le i,j \le n}$ is of full rank for all sufficiently small, positive values of the discretisation parameters.
\end{assumption}

\begin{theorem}\label{thm:HMod}
	Consider a truncation index $N\in \N$ such that \Cref{ass:SNonDeg} holds.
	There exists a Hamiltonian structure $(\mathrm{Jet}_1(Q) ,\omega,\mathcal H)$ such that Hamilton's equations recover the reduced equation \eqref{eq:2ndOrderODE} up to $\mathcal O(h^{N+1})$ terms. The symplectic structure $\omega$ and the Hamiltonian $\mathcal H$ can be chosen to be $\mathcal O(h)$-close to $\omega^0$ and $\mathcal H^0$, respectively.
\end{theorem}

\begin{proof}
	We use Ostrogradsky's Hamiltonian description of high-order Lagrangian systems (see, for instance, \cite{RASHID1996,Pons1989}).
	Let $M = \max\{M_1,\ldots,M_N\}$ denote the order of the highest derivative in the power series \eqref{eq:VarPrinc} truncated with $\mathcal{O}(h^{N+1})$ error. Denote the $(2M-1)$-jet space over $Q$ by $\mathrm{Jet}_{2M-1}(Q)$ with coordinates denoted by $\{\phi\}_{(2M-1)\mathrm{jet}}=(\phi,\dot \phi,\phi^{(2)},\ldots,\phi^{(2M-1)})$, where each $\phi^{(j)}$ is $\R^n$ valued. Let $Q^M$ denote the product $Q^M = Q \times Q \times \ldots \times Q$.
	We equip the cotangent bundle $\pi \colon T^\ast Q^M \to Q^M$ with Darboux coordinates $(q,p)=(q^1,\ldots,q^M,p_1,\ldots,p_M)$ such that $q^i = \phi^{(i-1)} \circ \pi$ and the symplectic structure is given by%
	\[ \Omega = \sum_{i=1}^M \sum_{j=1}^n \d q_j^i \wedge \d p^j_i.\]
	Consider the map $\overline \Psi \colon \mathrm{Jet}_{2M-1}(Q) \to T^\ast Q^M$ with
	$\{\phi\}_{(2M-1)\mathrm{jet}} \mapsto (q,p)$, where
	\begin{equation}\label{eq:defPsiOver}
		\begin{split}
			q^1 &= \phi\\[0.7em]
			&\vdots\\[0.7em]
			q^i &= \phi^{(i-1)} \\[0.7em]
			&\vdots\\[0.7em]
			q^{M} &= \phi^{(M-1)}
		\end{split}
		\qquad
		\begin{split}
			p^j_1 &= \frac{\p \mathcal{L}^{[N]} }{\p \dot \phi_j} - \frac{\d}{\d \xi} \frac{\p \mathcal{L}^{[N]}}{\p \ddot\phi_j} \ldots + (-1)^{M-1} \frac{\d^{M-1}}{\d \xi^{M-1}}\frac{\p \mathcal{L}^{[N]}}{\p \phi_j^{(M)}}\\
			&\vdots\\
			p_i^j &= \sum_{k=0}^{M-i} (-1)^k \frac{\d^{k}}{\d \xi^{k}}\frac{\p \mathcal{L}^{[N]}}{\p (\phi_j)^{(k+i)}}\\
			&\vdots\\
			p_M^j &= \frac{\p \mathcal{L}^{[N]}}{\p \phi_j^{(M)} }.
		\end{split}
	\end{equation}
	Here $\mathcal{L}^{[N]}$ denotes the truncation of the series $\mathcal{L}$ after $\mathcal{O}(h^N)$ terms. (Although not reflected in the notation, the map $\overline{\Psi}$ depends on the truncation index $N$ as well.)
		\Cref{ass:SNonDeg} guarantees that $\overline \Psi$ has a local inverse: given $(q,p) \in T^\ast Q^M$, $\{\phi\}_{(M-1)\mathrm{jet}}$ can immediately be recovered from $q$ by \eqref{eq:defPsiOver}. Then the relation
		\[
		p_M =  \nabla_{\phi^{(M)}} \mathcal{L}^{[N]}
		\]
		can locally be solved for $\phi^{(M)}$ in terms of $(q,p)$ by the implicit function theorem.
		In the expression
		\begin{equation}\label{eq:piVec}
			p_i
			= \sum_{k=0}^{M-i} (-1)^k
			\frac{\d^{k}}{\d \xi^{k}} \nabla_{\phi^{(k+i)}}  \mathcal{L}^{[N]}
		\end{equation}
		the $k$th summand depends on $\{ \phi\}_{(M+k)\mathrm{jet}}$ such that $p_i$ depends on $\{ \phi\}_{(2M-i)\mathrm{jet}}$. The coefficient matrix of $\phi^{(2M-i)}$ is given by $\mathrm{Hess}_{\phi^{(M)}}(\mathcal{L}^{[N]})$, which is invertible by \Cref{ass:SNonDeg}. The above equations \eqref{eq:piVec} with $i=M-1,\ldots,1$ can, therefore, be solved sequentially for $\phi^{(M+1)},\ldots,\phi^{(2M-1)}$.
		
		Let us denote the local expression for $\phi^{(M)}$ by $\tilde{\phi}^{(M)}(q,p)$.
		Consider the Hamiltonian system $(T^\ast Q^M, \Omega,H^{[N]})$ with Hamiltonian
		\begin{equation}\label{eq:HamH}
			\begin{split}
				H^{[N]}(q,p) &= \sum_{i=1}^M \langle p_i , \dot{q}^i \rangle_{\R^n} - \mathcal{L}^{[N]}\left(q^1,\ldots,q^M,\tilde{\phi}^{(M)}(q,p)\right)\\
				&= \sum_{i=1}^{M-1} \langle p_i , q^{i+1} \rangle_{\R^n} + \langle p_M,\tilde{\phi}^{(M)}(q,p)\rangle_{\R^n}  - \mathcal{L}^{[N]}\left(q^1,\ldots,q^M,\tilde{\phi}^{(M)}(q,p)\right).
			\end{split}
		\end{equation}
		On the right hand side of \eqref{eq:HamH} terms of order $\mathcal{O}(h^{N+1})$ are truncated.
		The equations of motions of a pullback $(\mathrm{Jet}_{2M-1}(Q), \overline{\Psi}^\ast\Omega,H^{[N]} \circ \overline{\Psi})$ of the Hamiltonian system $(T^\ast Q^M, \Omega,H^{[N]})$ via $\overline{\Psi}$ yield a jet extension of \eqref{eq:truncatedODE}.
		
		
		Consider coordinates $(\underline{\phi},\underline{\dot\phi})$ on the 1-jet bundle $\underline\pi\colon \mathrm{Jet}_1(Q) \to Q$ such that $\underline{\phi} = \phi \circ \underline \pi$, where $\underline\pi$ is the jet bundle projection. We define the inclusion
		
		\[\Psi \colon \mathrm{Jet}_1(Q) \to \mathrm{Jet}_{2M-1}(Q) \quad \text{by} \quad
		\phi^{(j)} \circ \Psi = \begin{cases}
			\underline{\phi} \quad &\text{if } j=0\\
			\underline{\dot \phi} \quad &\text{if } j=1\\
			g_j(\underline{\phi},\underline{\dot\phi}) \quad &\text{otherwise},
		\end{cases}
		\]
		where $g_j(\underline{\phi},\underline{\dot\phi})$ is the substitution of the $j$-th derivative considered in the order reduction process in \eqref{eq:truncatedODE} and \eqref{eq:2ndOrderODE} truncated after $\mathcal{O}(h^N)$ terms.
		Consider the pullback of the Hamiltonian system $(\mathrm{Jet}_{2M-1}(Q), \overline{\Psi}^\ast\Omega,H^{[N]} \circ \overline{\Psi})$ via $\Psi$ to $(\mathrm{Jet}_1(Q),\omega^{[N]},\mathcal{H}^{[N]})$ with $\omega^{[N]} =  \Psi^\ast \overline{\Psi}^\ast\Omega$ and $\mathcal{H}^{[N]} = H^{[N]} \circ \overline{\Psi} \circ \Psi$.
		For $h$ close to $0$ the 2-form $\omega^{[N]}$ is non-degenerate due to the structure of $\Psi$. As pullback and an application of $\d$ commute, $\omega^{[N]}$ is a symplectic form and the range $\rg(\Psi)$ of $\Psi$ is a symplectic submanifold of $(\mathrm{Jet}_{2M-1}(Q), \overline{\Psi}^\ast\Omega)$.
		Let $\pi \colon (\mathrm{Jet}_{2M-1}(Q) \to \mathrm{Jet}_{1}(Q)$ denote the projection to 1-jets and let $\pi' = \psi \circ \pi \colon \mathrm{Jet}_{2M-1}(Q) \to \rg(\Psi)$.
		Let $\overline X^{[N]}$ denote the Hamiltonian vector field of the system $(\mathrm{Jet}_{2M-1}(Q), \overline{\Psi}^\ast\Omega,H^{[N]} \circ \overline{\Psi})$ and $\mathcal X^{[N]}$ of the system $(\mathrm{Jet}_1(Q),\omega^{[N]},\mathcal{H}^{[N]})$.
		The restriction of $\pi'_\ast \overline X^{[N]}$ to the image of $\Psi$ is $\Psi$-related to the vector field $\mathcal X^{[N]}$, i.e.\
		\[
		(\Psi|_{\rg(\Psi)}^{-1})_\ast (\pi'_\ast \overline X^{[N]}\circ \psi)  = \mathcal X^{[N]}.
		\]
		Here $\Psi|_{\rg(\Psi)}^{-1}$ denotes the inverse of $\Psi$ considered as the diffeomorphism $\Psi \colon \mathrm{Jet}_1(Q) \to \rg(\Psi)$.
		As the flow of $\overline{X}^{[N]}$ leaves the image of $\Psi$ invariant up to order $\mathcal{O}(h^{N+1})$-terms, it follows that the flow of $\overline{X}^{[N]}$ and $\mathcal X^{[N]}$ are $\Psi$-related up to terms of order $\mathcal{O}(h^{N+1})$ (as explicitly argued in \cite[Prop.5.1]{BEAMulti}).
	\end{proof}

	\begin{corollary}\label{thm:LMod2ndOrder}
		Consider a truncation index $N\in \N$ such that \Cref{ass:SNonDeg} holds. Let $(\mathrm{Jet}_1(Q) ,\omega,\mathcal H)$ be the Hamiltonian structure provided by  \cref{thm:HMod}. If there exist Darboux coordinates $(q,p)$ on $(\mathrm{Jet}_1(Q),\omega)$ that are $\mathcal{O}(h)$-close to $(\mathfrak q, \mathfrak p)$ such that $\left(\frac{\p^2 \mathcal H}{\p{p}_i \p{p}_j }\right)_{1\le i,j\le n}$ is of full rank, then there exists a first-order Lagrangian $\tilde {\mathcal L}$ given as
		\begin{equation}\label{eq:AssFormL}
			\begin{split}
				\tilde {\mathcal L}(q,\dot q) &= \tilde {\mathcal L}^0(q,\dot q)+ h \tilde {\mathcal L}^1(q,\dot q) + h^2 \tilde {\mathcal L}^2(q,\dot q)+\ldots + h^N \tilde {\mathcal L}^N(q,\dot q)\\
				&= \sum_{i=0}^N h^i \tilde {\mathcal L}^i(q,\dot q)
			\end{split}
		\end{equation}
		such that $\tilde {\mathcal L}^0=\mathcal L^0$ and the Euler--Lagrange equations $\mathcal E_2 \tilde {\mathcal L} =0$ recover \eqref{eq:2ndOrderODE} up to $\mathcal{O}(h^{N+1})$-terms.
	\end{corollary}
	
	\begin{proof}
		The Lagrangian $\tilde {\mathcal L}$ is obtained as the Legendre transform of $\mathcal H$, i.e.\
		\[\tilde {\mathcal L} =\sum_{i=1}^n{\dot q}^i {p}_i - {\mathcal H},\]
		where all quantities are expressed in $({q},{\dot q})$. The zeroth order term in $h$ coincides with $\mathcal L^0$ by construction.
	\end{proof}

	\begin{remark}\label{rem:coordFree1stOrder}
		The non-degeneracy assumption on $\mathcal H$ and the choice of Darboux coordinates can be dropped in a coordinate-free description of the motion. Here {\em first-order principle} refers to a variational principle with a 1-form defined on a 1-jet space.
		As the symplectic form $\omega$ is closed, it has a local primitive $\lambda$ which we can chose to be $\mathcal O(h)$ close to $\lambda^0 = -\sum_{i=1}^n\mathfrak p_i \d \mathfrak q^i$.
		The Lagrangian density $\mathfrak L$ is given as the 1-form $\lambda  -\mathcal{H} \d \xi$ on $\mathrm{Jet}_1(Q) \times \R$. The corresponding variational principle $S(\gamma) = \int_a^b \gamma^\ast \mathfrak L$, where $\gamma\colon [a,b]\to \mathrm{Jet}_1(Q)$ is the prolongation of a $Q$-valued curve, is of first order.
	\end{remark}
	
	\begin{remark}\label{rem:NoBundle}
		If we transform the coordinates $({q},{p})$ from \cref{thm:LMod2ndOrder} back to the dynamical coordinates $(\phi,\dot \phi)$ then ${q}$ will typically depend on $\dot \phi$ as well. Therefore, an expression of the Lagrangian $\tilde {\mathcal L}$ in the original variables can involve second derivatives of $\phi$. This is because the map $\Psi \colon \mathrm{Jet}_1(Q) \to \mathrm{Jet}_{2M-1}(Q)$ does not respect the bundle structure of the jet spaces. Therefore, the distribution $\mathcal D$ spanned by the vector fields $\frac{\p}{\p \dot \phi_j}$ might not be Lagrangian for the structure $\omega$ and there exists no primitive $\lambda$ with kernel $\mathcal D$. Therefore, an expression of $\lambda$ in the frame $\d \phi_j,\d \dot{\phi}_j$ must involve $\d \dot{\phi}_j$ components. An expression of the variational principle $S(\gamma) = \int_a^b \gamma^\ast \mathfrak L$ from \cref{rem:coordFree1stOrder} in coordinates will involve second derivatives of the curve $\gamma$.
		If $\mathcal D$, however, happens to be Lagrangian, then $\omega$ admits a primitive $\lambda$ with kernel $\mathcal D$ \cite[Cor.\ 15.7]{Libermann1987}. The variational principle constructed in \cref{rem:coordFree1stOrder} to $\lambda$ is then expressible in the coordinates $(\phi,\dot \phi)$.
	\end{remark}
	
	\begin{remark}
		To obtain $L(\phi,\dot \phi)$ from $\mathcal{L}(\{\phi\}_{\infty\mathrm{jet}})$ it is in general {\em not} possible to simply substitute the higher order derivatives in the Lagrangian $\mathcal{L}(\{\phi\}_{\infty\mathrm{jet}})$ using \eqref{eq:truncatedODE} and its derivatives.
		Indeed, the substitution on the Lagrangian side only works if $\mathcal L$ has the form of a meshed Lagrangian (see \cite{Vermeeren2017}), which cannot be assumed in this context. We have seen, however, that the substitution can safely be done on the Hamiltonian side.
	\end{remark}

	\begin{remark}
		As conserved quantities of the original system \eqref{eq:formalODE} are passed on to the reduced system by \cref{prop:presConsQ}, (infinitesimal) symmetries in the sense of Noether's theorem are passed on to the Hamiltonian system $(\mathrm{Jet}_1(Q) ,\omega,\mathcal H)$ as well as to the modified Lagrangian density.
	\end{remark}
	
	\begin{proof}[Proof of \cref{thm:ConsPalais}]
		An application of \cref{thm:LMod2ndOrder} to a situation where the principle \eqref{eq:VarPrinc} arises as the series expansion of a symmetrised discrete Lagrangian for a symmetry group of codimension 1 yields the statement about variational integrators formulated in \cref{thm:ConsPalais}.
	\end{proof}
	
	%

	\section{Rotating travelling waves in the discretised nonlinear wave equation}\label{sec:TravWave}
	
	Let us consider travelling waves with constant phase rotation in the nonlinear wave equation \eqref{eq:NonLinearWavePDE}.
	Employing methods of the proof of \cref{thm:HMod}, we will compute a modified Hamiltonian system defined on a phase space of minimal dimension which governs the rotating travelling waves in the discretised equations for the five-point stencil.
	The modified Hamiltonian system corresponds to a first-order Lagrangian principle.
	For special cases such as no rotation, zero wave speed, or special choices of the steps sizes,
	the modified Lagrangian can be expressed in the original variables, i.e.\ the same variables as in the continuous setting.
	
	\subsection{Continuous setting}
	
	The Euler--Lagrange equation of the action
	
	\begin{equation}
		\label{eq:FunctionalSForRotEx}
		S(u) = \int \left(\frac 12 \left( \langle u_t , u_t \rangle - \langle u_x , u_x \rangle\right) + W(u) \right)\,\d t\, \d x
	\end{equation}
	is the nonlinear wave equation \eqref{eq:NonLinearWavePDE}, i.e.\
	
	\begin{equation}\label{eq:PDERotExample}
		u_{tt} - u_{xx} -\nabla W(u) = 0.
	\end{equation}
	
	In the following, we restrict ourselves to $W(u) = \frac 12 V( \langle u , u \rangle)$, $u \colon \R^2 \to \R^2$, and $V\colon\R\to\R$ analytic in order to analyse rotating travelling wave solutions. The considerations hold in the more general setting of any $W$ and $u \colon \R^2\to \R^d$ in the non-rotating case, and for any symmetries of the form $\partial_x - c\partial_t + g(u)_i\partial_{u_i}$.

	\begin{lemma}
		Solutions of \eqref{eq:PDERotExample} of the form
		\begin{equation}\label{eq:AnsSymSols}
			u(t,x) = R(t)\phi(x-ct)
		\end{equation}
		with $\phi \colon \R \to \R^2$ and
		\begin{equation}
			\label{eq:defJR}
			J = \begin{pmatrix}0&1\\ -1 & 0 \end{pmatrix},
			\qquad
			R(t) = \exp(t \alpha J) = \begin{pmatrix}
				\phantom{-}\cos(\alpha t) & \sin(\alpha t)\\
				-\sin(\alpha t) & \cos(\alpha t)\\
			\end{pmatrix}
		\end{equation}
		solve the ODE
		\begin{equation}\label{eq:SymmetrisedRotInv}
			(\alpha^2 + V'(\langle \phi(\xi), \phi(\xi) \rangle)) \phi(\xi) + 2 c \alpha J \dot{\phi}(\xi) - (c^2-1) \ddot{\phi}(\xi)=0.
		\end{equation}
		On the other hand, solutions to \eqref{eq:SymmetrisedRotInv} give rise to solutions $u(t,x) = R(t)\phi(x-ct)$ of \eqref{eq:PDERotExample}.
	\end{lemma}
	
	\begin{proof}
		A substitution of $u(t,x) = R(t)\phi(x-ct)$ into \eqref{eq:PDERotExample} yields
		\[
		(c^2-1) R(t) \ddot{\phi}(\xi) - (\alpha^2 + V'(\langle \phi(\xi), \phi(\xi) \rangle))R(t) \phi(\xi) - 2 \alpha c J R(t) \dot{\phi}(\xi) =0,
		\]
		where $\xi = x-ct$. Since $J$ commutes with $R(t)$ and $R(t)$ is invertible, this is equivalent to \eqref{eq:SymmetrisedRotInv}.
	\end{proof}

	\begin{figure}
		\includegraphics[width=0.4\textwidth]{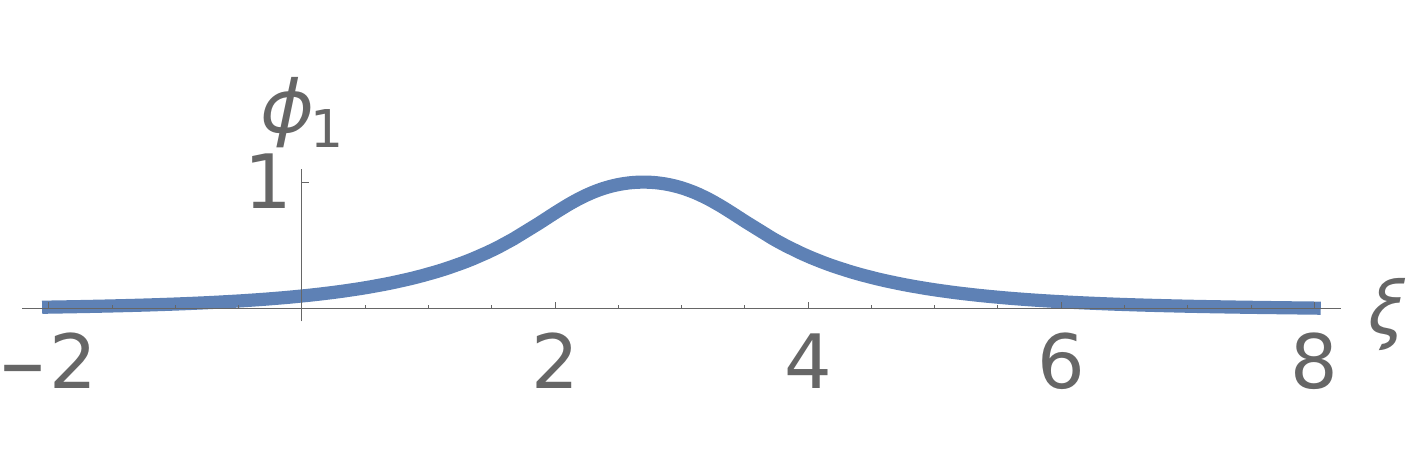}
		\includegraphics[width=0.4\textwidth]{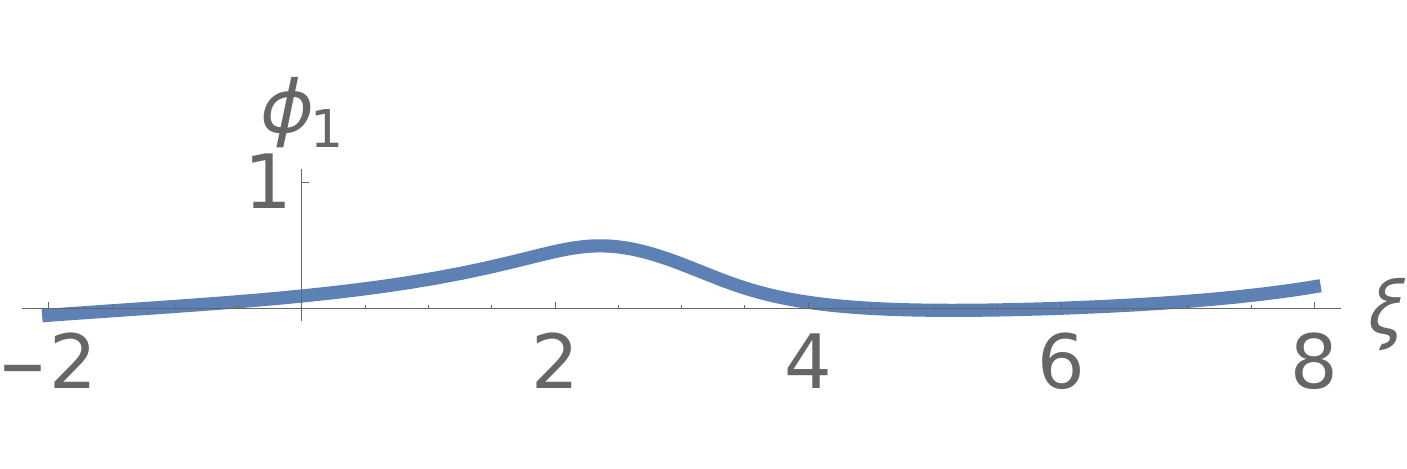}\\
		\includegraphics[width=0.4\textwidth]{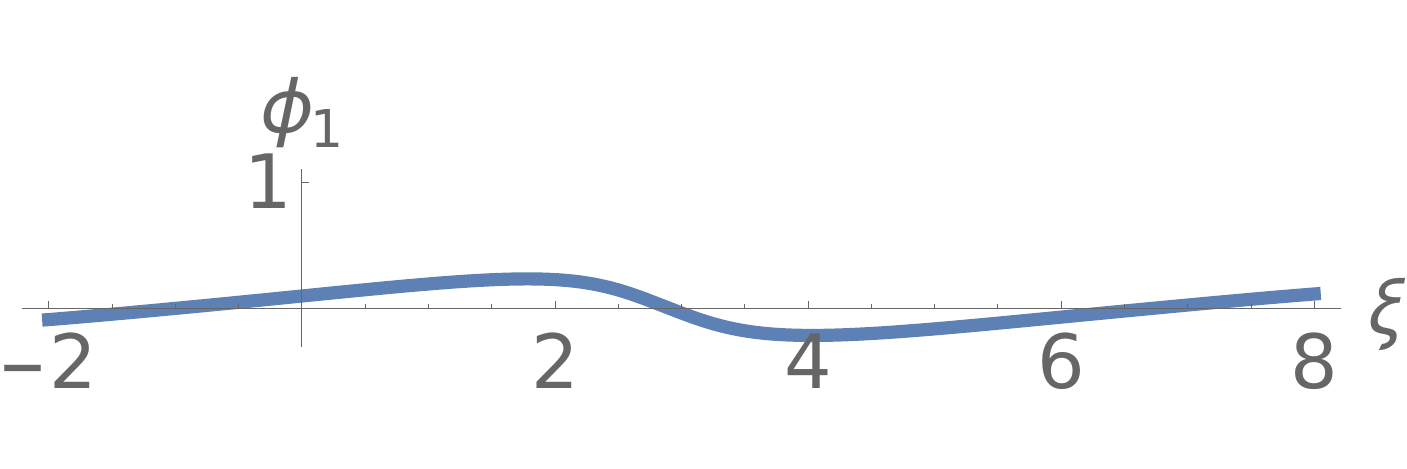}
		\includegraphics[width=0.4\textwidth]{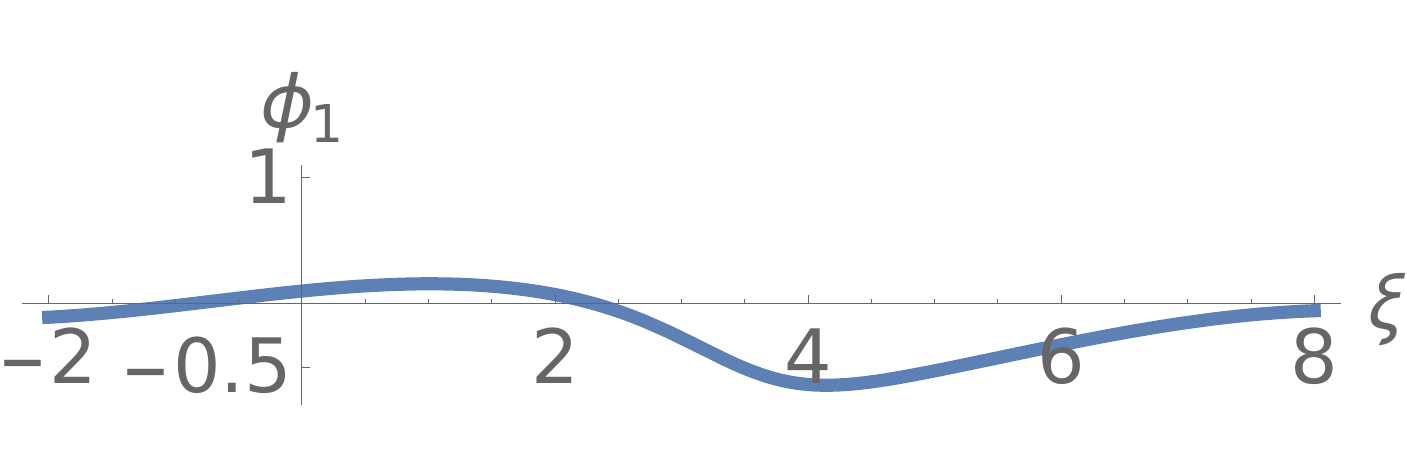}
		\caption[Dynamics of amplitude variable]{Dynamics of the amplitude variable $\phi_1(\xi)$ for $\alpha \in \{0,0.3,0.5,0.7\}$ for $V(a) = -\exp(-(a-1)^2)$ and the wave speed $c=0.5$. (Initial condition $\phi_1(0)=\phi_2(0)=\dot\phi_1(0)=\dot\phi_2(0)=0.1$)}\label{fig:dynamicsPhi1}
	\end{figure}
	
	\begin{figure}
		\begin{center}
			\includegraphics[width=0.2\textwidth]{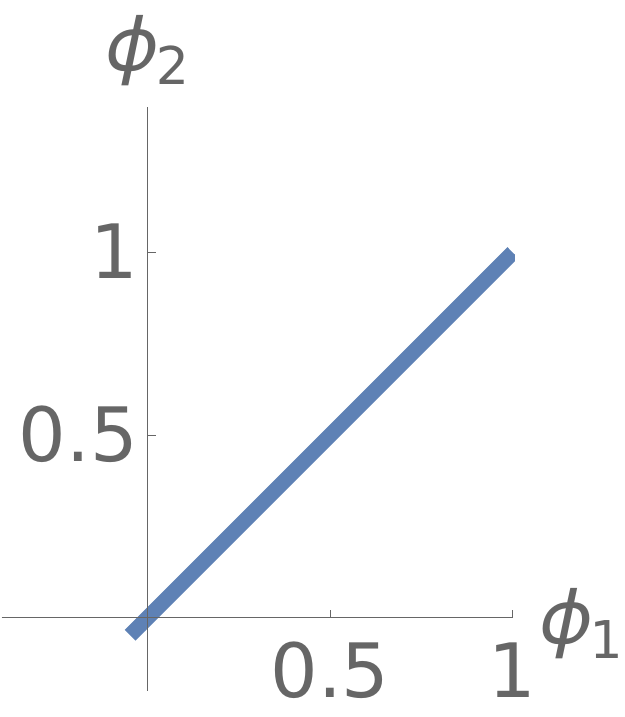}
			\includegraphics[width=0.2\textwidth]{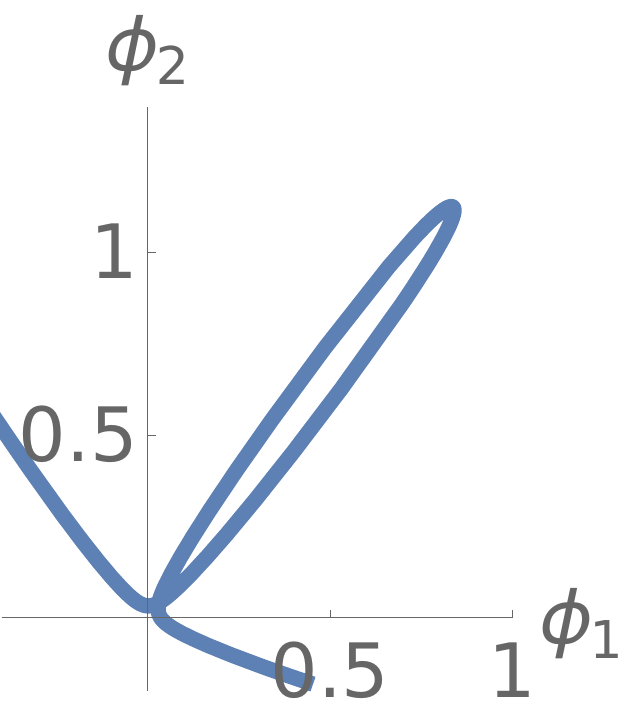}
			\includegraphics[width=0.2\textwidth]{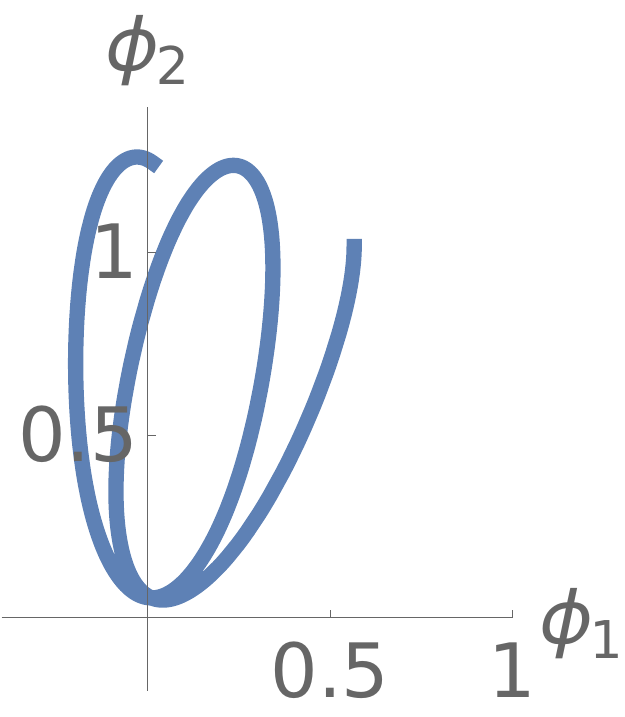}
			\includegraphics[width=0.2\textwidth]{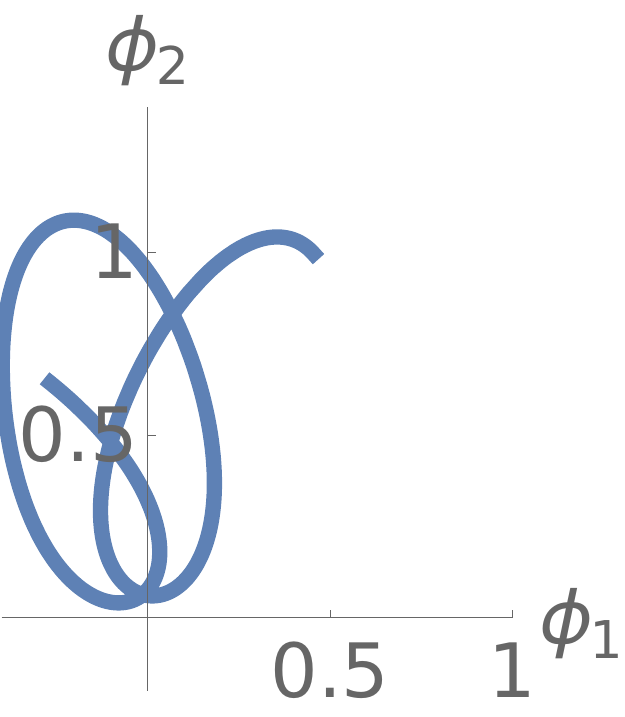}
		\end{center}
		\caption[Phase portrait amplitude variable]{Phase portrait of the amplitude variables $\phi_1(\xi)$, $\phi_2(\xi)$ for $\alpha \in \{0,0.1,0.6\}$, $V(a) = -\exp(-(a-1)^2)$, the wave speed $c=0.5$ and $\xi \in [-5,10]$. (Initial condition $\phi_1(0)=\phi_2(0)=\dot\phi_1(0)=\dot\phi_2(0)=0.1$)}\label{fig:PhasePhi1Phi2}
	\end{figure}
	
	The dynamics of $\phi_1(\xi)$ for different values of $\alpha$ and $V(a)=-\exp(-(a-1)^2)$ are displayed in \cref{fig:dynamicsPhi1}. \Cref{fig:PhasePhi1Phi2} shows phase plots of $\phi_1(\xi)$ and $\phi_2(\xi)$.

	\begin{remark}\label{rem:criticalwavespeed}
		There is a critical wave speed $c \in \{-1,1\}$ for which the ODE \eqref{eq:SymmetrisedRotInv} is of order 1 if the rotational speed $\alpha \not = 0$ and algebraic if $\alpha =0$. It is not be governed by a non-degenerate, autonomous Lagrangian and cannot be analysed within the autonomous Hamiltonian framework. The case $c^2=1$ is excluded in the following.
	\end{remark}
	
	\begin{lemma}\label{lem:Ssym}
		The system of ODEs \eqref{eq:SymmetrisedRotInv} are the Euler-Lagrange equations to the action functional
		\[
		S^{\mathrm{sym}}(\phi) = \int L^0(\phi,\dot \phi) \, \d \xi
		\]
		with
		\begin{equation}\label{eq:L0Def}
			L^0(\phi,\dot \phi) =\frac{1}{2} \left( \alpha^2 \langle \phi,\phi\rangle - 2 \alpha c \langle J \phi, \dot{\phi}\rangle + (c^2-1) \langle \dot{\phi},\dot{\phi}\rangle + V( \langle \phi,\phi\rangle ) \right).
		\end{equation}
	\end{lemma}
	%
	
	\begin{remark}
		
		Indeed, restricting $S$ to symmetric functions of the form $(t,\xi) \mapsto R(t) \phi(\xi)$ with $\xi = x-ct$ yields the functional $S^\mathrm{sym}$ from \cref{lem:Ssym}. This shows that Palais' principle of symmetric criticality \cite{palais1979} is valid in this example, i.e.\ the critical points of $S$ which are symmetric coincide with the points which are symmetric and critical with respect to symmetric variations. In other words, if $u$ is symmetric and $\D S(u)(v)=0$ for all symmetric test functions $v$ then $\D S(u)(v)=0$ for all test functions $v$. Here, we assume that $S$ can be defined on a Banach space\footnote{The exact set-up may depend on $V$.} and $\D$ is the Fr\'echet derivative.
		The validity of the principle of symmetric criticality can then be concluded directly from \cite{palais1979} using the action of the compact Lie group $(\R / \frac{2\pi}{\alpha}\Z)$ on the domain of definition of $S$ given by $(s \cdot u)(t,\xi) = R(-s)u(t+s,\xi)$, where $u$ is expressed in the coordinates $(t,\xi)$ with $\xi=x-ct$.
	\end{remark}

	\begin{remark}
		The ODE \eqref{eq:SymmetrisedRotInv} admits a Hamiltonian formulation on $\R^4$ equipped with the symplectic structure
		\[
		\d {\mathfrak q}^1 \wedge \d \mathfrak p_1 + \d \mathfrak q^2 \wedge \d \mathfrak p_2.
		\]
		The super- or subscripts of $\mathfrak q$ and $\mathfrak p$ denote components of $\mathfrak q$ or $\mathfrak p$, respectively.
		The Hamiltonian is obtained as the Legendre transformation of the Lagrangian function $L^0$ defined in \eqref{eq:L0Def} and expressed in the canonical coordinates
		\begin{equation}\label{eq:coordsQP}
			\mathfrak q = \phi, \qquad \mathfrak p = \nabla_{\dot \phi}L^0 = (c^2-1) \dot  \phi - c \alpha J  \phi.
		\end{equation}
		The Hamiltonian is given as
		\[
		\mathfrak{H}(\mathfrak q, \mathfrak p) = \frac 1 {2(c^2-1)} (\| \mathfrak p\|^2+2 c \alpha \langle \mathfrak p,J \mathfrak q\rangle + \alpha^2 \|\mathfrak q\|^2-(c^2-1)V(\| \mathfrak q\|^2)).
		\]
	\end{remark}

	\begin{remark}\label{rem:RotInvExact}
		The 1-form $L^0(\phi,\dot \phi) \d \xi$ is invariant under the prolongation of the Lie group action of $S^1 \cong \R / (2 \pi \Z)$ on $\R^2$ defined by $\theta \cdot (\xi,\phi) = (\xi, \exp(\theta J^\top) \phi)$, where $\exp$ denotes the matrix exponential as in \eqref{eq:defJR}.
		By Noether's First Theorem \cite[\S 7.2]{mansfield2010} the quantity
		\[
		I^{\mathrm{dyn}}_{\mathrm{rot}}(\phi,\dot \phi) = \langle \nabla_{\dot \phi}L^0,J^\top \phi \rangle= \alpha c \|\phi\|^2 + (c^2-1) \langle J \dot \phi, \phi \rangle
		\]
		is conserved along solutions of \eqref{eq:SymmetrisedRotInv}. In the canonical coordinates $\mathfrak q$, $\mathfrak p$ from \eqref{eq:coordsQP} the quantity is given as
		\[
		I_{\mathrm{rot}}(\mathfrak q, \mathfrak p) = \langle J \mathfrak p, \mathfrak q \rangle.
		\]
\begin{figure}
		\begin{center}
			\includegraphics[width=0.4\textwidth]{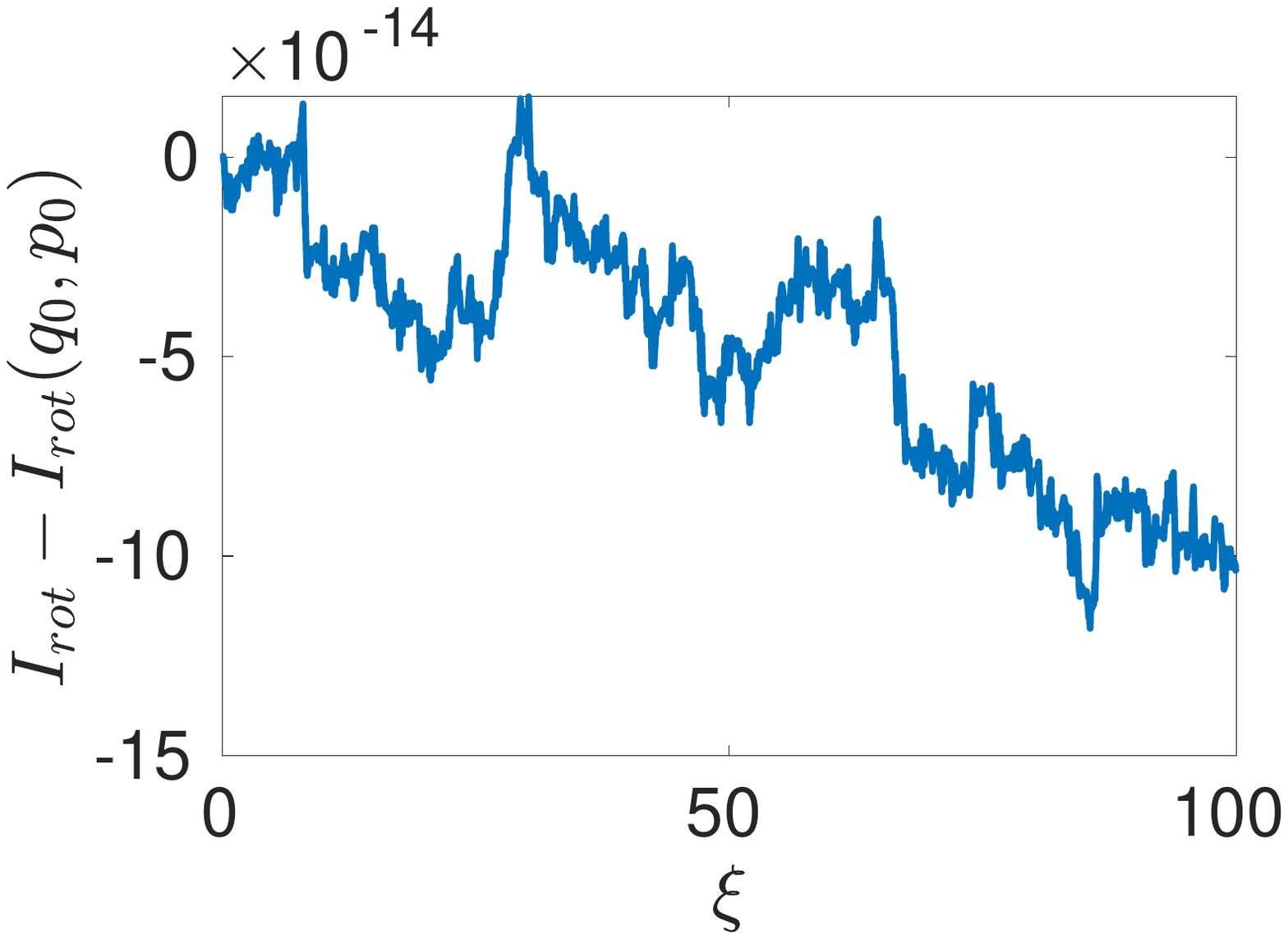}
		\end{center}
		\caption[Numerical preservation of rotational quantity]{
			Evaluation of the conserved quantity $I_{\mathrm{rot}}$ (see \cref{rem:RotInvExact}) along a numerically computed trajectory shows round-off errors only (vertical axis is scaled by $10^{-14}$). Here $V(a)= -\frac 12 a -a^2$, $\alpha = -1$, $c=2$. The integrator is the symplectic midpoint rule. Implicit equations are solved using fixed point iterations.
		}\label{fig:IRotExact}
	\end{figure}
		Also see \cref{fig:IRotExact}. The quantities $\mathfrak H$ and $I_{\mathrm{rot}}$ are functionally independent on a superset of the dense open subset $\{(\mathfrak q, \mathfrak p) | \langle \mathfrak p,\mathfrak q\rangle \not =0 \}$ of the phase space. As they also Poisson commute, the considered system \eqref{eq:SymmetrisedRotInv} is Liouville completely integrable.
	\end{remark}

	\subsection{Five-point stencil discretisation and modified equation}
	
	The five-point stencil discretisation of \eqref{eq:PDERotExample} with respect to a mesh $\{(i \Delta t, j \Delta x)\}_{(i,j)\in \Z^2}$ and $u_{i,j}$ corresponding to the value of a function $u$ at the meshpoint $(i \Delta t, j \Delta x)$ is given as
	
	\begin{equation}
		\label{eq:PDERotExample5PtStencilClassic}
		0=\frac{1}{{\Delta t}^2} \left( u_{i-1,j}  - 2 u_{i,j} + u_{i+1,j} \right)
		- \frac{1}{{\Delta x}^2} \left( u_{i,j-1} - 2 u_{i,j} + u_{i,j+1} \right)
		-V'(\langle u_{i,j}, u_{i,j}\rangle) u_{i,j}.
	\end{equation}
	
	The scheme is multisymplectic. It arises via a discretisation of the continuous action $S$ as the following lemma shows.
	\begin{lemma}
		A discrete solution $u_\Delta = (u_{i,j})_{i,j \in \Z}$ satisfies \eqref{eq:PDERotExample5PtStencilClassic} if and only if for all $K \in \N$ it extremises
		\begin{equation*}
			S^K_{\Delta}(u) = \frac{1}{2}\sum_{i,j=-K}^K \frac{\|u_{i-1,j}-u_{i,j}\|^2}{\Delta t^2} - \frac{\|u_{i,j-1}-u_{i,j}\|^2}{\Delta x^2}- V(\|u_{i,j}\|^2 )
		\end{equation*}
		on all interior grid-points, i.e.\ $\nabla_{(u_{i,j})_{-K+1\le i,j \le K-1}} S^K_\Delta (u) = 0$.
	\end{lemma}
	
	As discussed in the introduction, we pass to the functional equation
	\begin{equation}
		\label{eq:PDERotExample5PtStencil}
		\begin{aligned}
			0&=\frac{1}{{\Delta t}^2} \left( u(t-\Delta t, x) - 2 u(t,x) + u(t+\Delta t,x) \right) \\
			&- \frac{1}{{\Delta x}^2} \left( u(t, x - \Delta x) - 2 u(t,x) + u(t,x+\Delta x) \right) \\
			&-V'(\langle u(t,x), u(t,x)\rangle) u(t,x)
		\end{aligned}
	\end{equation}
	with $(t,x)\in\R^2$, whose solutions coincide with \eqref{eq:PDERotExample5PtStencilClassic} on the grid.
	
	The ansatz for a symmetric solution from \eqref{eq:AnsSymSols}, i.e.\ $u(t,x) = R(t) \phi(x-ct)$ with $\xi = x-ct$, leads to the following functional equation for $\phi$
	\begin{equation}\label{eq:5PtStencilSymmetric}
		\begin{aligned}
			0&=\frac{1}{h^2{\Delta t}^2} \left( R(-h\Delta t) \phi(\xi + c h\Delta t) - 2 \phi(\xi) + R(h\Delta t) \phi(\xi - c h\Delta t) \right)\\
			&-\frac{1}{h^2{\Delta x}^2} \left( \phi(\xi + h\Delta x) - 2 \phi(\xi) + \phi(\xi - h\Delta x) \right)\\
			&-V'(\langle \phi(\xi), \phi(\xi)\rangle) \phi(\xi).
		\end{aligned}
	\end{equation}
	Here we have introduced the formal series variable $h$ to the same power as the step sizes.
	A series expansion of \eqref{eq:5PtStencilSymmetric} around $h=0$ followed by solving for $\ddot{\phi}$ in terms of $\phi$, $\dot \phi$ and higher order terms yields a formal power series of the form
	\begin{equation}
		\label{eq:Taylor5PtStencilSymmetric}
		\begin{aligned}
			\ddot{\phi}(\xi) &= \frac{(\alpha^2 + V'(\langle \phi(\xi), \phi(\xi) \rangle)) \phi(\xi) + 2 c \alpha J \dot{\phi}(\xi)}{c^2-1}\\
			&+ h^2  g_2(\phi^{(4)}(\xi),\ldots,\dot{\phi}(\xi),\phi(\xi)) \\
			&+ h^4  g_4(\phi^{(6)}(\xi),\ldots,\dot{\phi(\xi)},\phi(\xi)) \\
			&+\ldots.
		\end{aligned}
	\end{equation}
	Recall that the critical wave speed $c^2=1$ is excluded from our discussion (see \cref{rem:criticalwavespeed}).
	Using \eqref{eq:Taylor5PtStencilSymmetric} to replace $\ddot{\phi}$ and all higher derivatives on the right hand side of \eqref{eq:Taylor5PtStencilSymmetric} makes second order derivatives only occur in $h^4$ and higher order terms. Repeating this process iteratively we can push derivatives of order greater than 2 to $\mathcal{O}(h^r)$  for arbitrary $r$.
	We obtain a formal series of the form
	\begin{equation}
		\label{eq:Phi2Elliminated}
		\ddot{\phi}(\xi) = \frac{(\alpha^2 + V'(\langle \phi(\xi), \phi(\xi) \rangle)) \phi(\xi) + 2 c \alpha J \dot{\phi}(\xi)}{c^2-1}
		+\sum_{j=1}^{\infty} h^{2j}  \hat g_{2j}(\dot{\phi}(\xi),\phi(\xi)).
	\end{equation}
	\begin{figure}
		\includegraphics[width=0.4\textwidth]{Phi1a0_bigLabels.pdf}
		\includegraphics[width=0.4\textwidth]{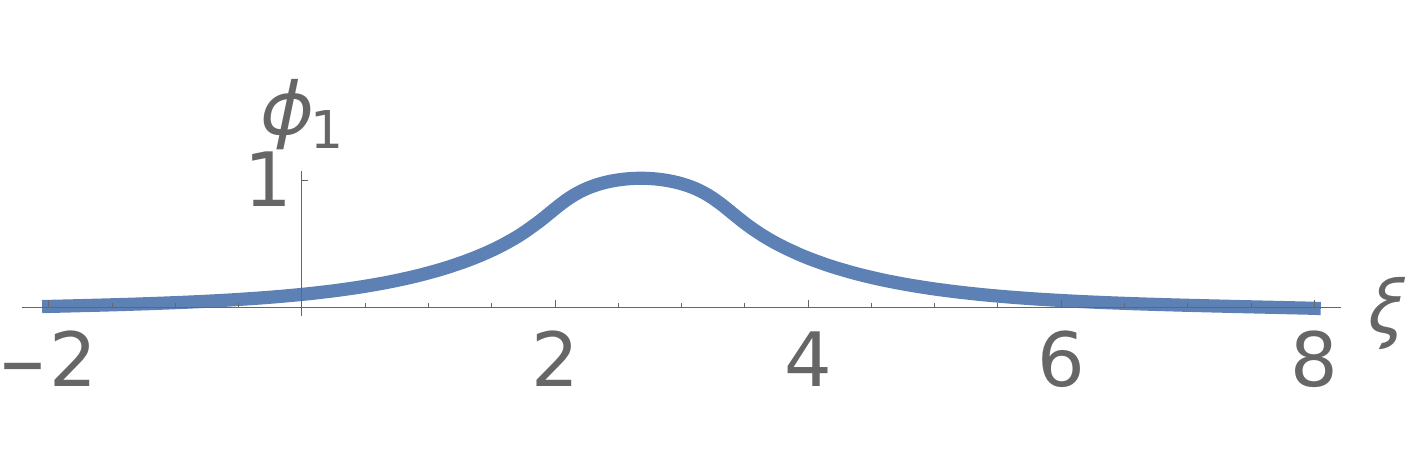}\\
		\includegraphics[width=0.4\textwidth]{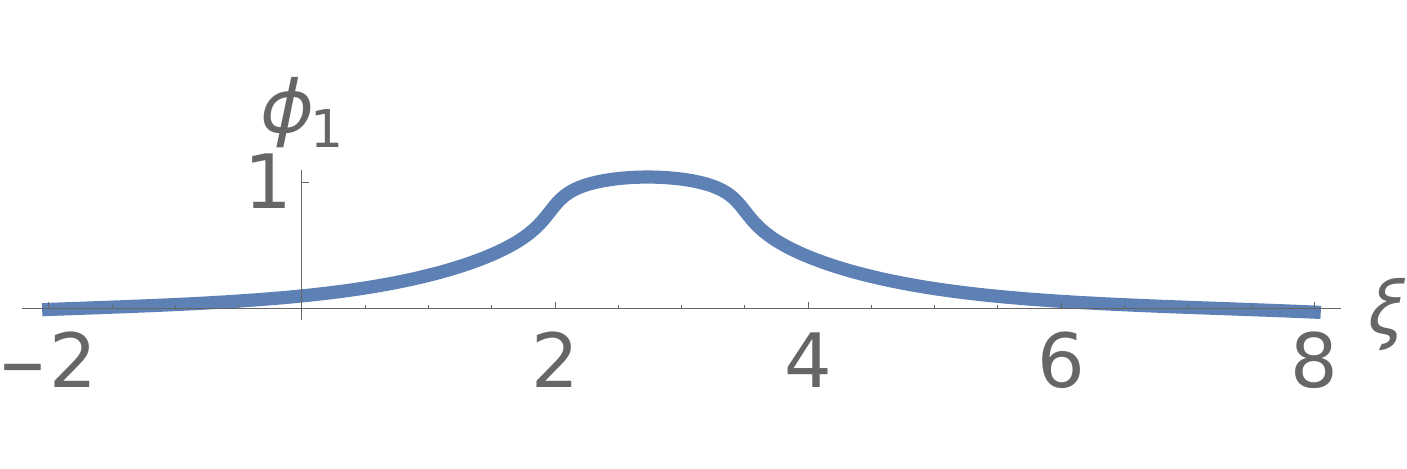}
		\includegraphics[width=0.4\textwidth]{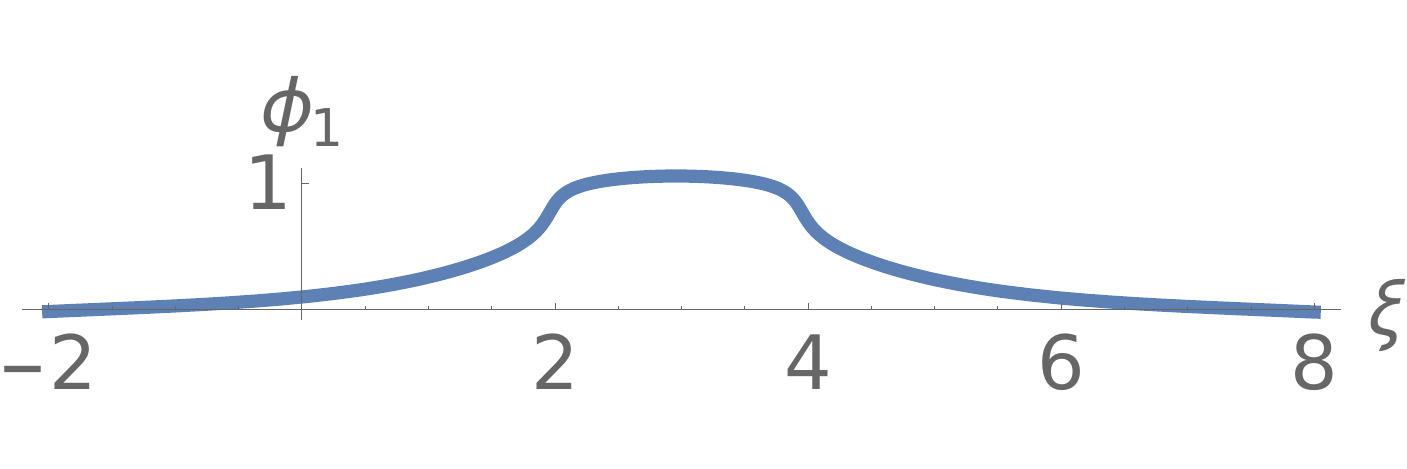}
		\caption[Dynamics modified equation]{Dynamics of the amplitude variable $\phi_1(\xi)$ for $\alpha=0$, $V(a) = -\exp(-(a-1)^2)$, $c=0.5$ and $\Delta x \in \{0,0.6,1,1.2\}$ for the modified equation truncated after $\mathcal O(h^3)$ terms.}\label{fig:dynamicsPhi1dx}
	\end{figure}
	The second order term $\hat g_{2}$ is reported in {\tt  Computational\_Results\_documented.pdf} in \cite{multisymplecticSoftware}. The dependence of the dynamics of $\phi_1(\xi)$ on the step size is illustrated in \cref{fig:dynamicsPhi1dx}.

	\begin{remark}\label{rem:symDiscr}
		The considered discretisation respects the rotation symmetry $\theta \cdot (\xi,\phi) = (\xi, \exp(\theta J^\top) \phi)$ introduced in \cref{rem:RotInvExact}. Therefore, \eqref{eq:5PtStencilSymmetric}, \eqref{eq:Taylor5PtStencilSymmetric}, \eqref{eq:Phi2Elliminated} are invariant under a prolongation of the action.
	\end{remark}
	
	\subsection{Computation of modified structures}\label{sec:structMod}

	We now follow the proof of \cref{thm:HMod}. To simplify notation, we neglect to include the order $N$ of truncation in $h$ when denoting $\mathcal{L}$, $H$, $\mathcal{H}$, and $\omega$. Notice that a computation with $\mathcal{L}$ truncated to higher order than $N$ recovers all terms of $\mathcal{L}$, $H$, $\mathcal{H}$, and $\omega$ up to order $N$.
	However, we require the truncation order to be consistent.
	In the following computational example, we truncate after order $h^5$-terms.
	
	\begin{itemize}
		\item
		We compute the Hamiltonian system $(T^\ast Q^M,\Omega,H)$ governing the high-order equation \eqref{eq:Taylor5PtStencilSymmetric} on a sufficiently large phase space with the canonical symplectic structure $\Omega$.
		\item
		Expressing the above Hamiltonian system in dynamical coordinates then corresponds to the computation of $(\mathrm{Jet}_{2M}(Q), \overline{\Psi}^\ast\Omega,H \circ \overline{\Psi})$.
		\item
		We will then substitute on-shell solutions of higher order derivatives of $\phi$ (just as we did to obtain \eqref{eq:Phi2Elliminated} from \eqref{eq:Taylor5PtStencilSymmetric}) into the coordinates of the 2-form $\overline{\Psi}^\ast\Omega$ and the Hamiltonian $H \circ \overline{\Psi}$ to obtain the reduced modified Hamiltonian system $(\mathrm{Jet}_1(Q),\omega,\mathcal H)$ which governs $\eqref{eq:Phi2Elliminated}$.
		\item
		While the reduced modified Hamiltonian system can be explicitly computed, the corresponding modified first-order Lagrangian will, in the most general case, not have a closed form. However, we compute the expressions for special cases.
	\end{itemize}

	A series expansion of the discrete Lagrangian
	\[
	L_\Delta = \frac{\|u(t-h\Delta t,x)-u(t,x)\|^2}{\Delta t^2} - \frac{\|u(t,x-h\Delta x)-u(t,x)\|^2}{\Delta x^2}-\frac 12 V(\|u(t,x)\|^2 )
	\]
	followed by the substitution $u(t,x) = R(t)\phi(x-ct)$ yields

	\allowdisplaybreaks
	\begin{align*}
		\begin{autobreak}
			\mathcal L
=
\frac{1}{2} (\alpha^2 \|\phi \|^2
+2 \alpha c (\phi _1 \dot \phi _2- \phi _2 \dot \phi _1)
+(c^2-1)\|\dot \phi\|^2
-V)
+\frac{h^2}{24} \Big(\alpha^4 (
-{\Delta t}^2) \phi _2^2
-\alpha^4 {\Delta t}^2 \phi _1^2
+2 \alpha c {\Delta t}^2 \phi _2 (2 \alpha^2 \dot \phi _1
+3 \alpha c \ddot \phi _2
-2 c^2 \phi _1^{(3)})
+2 \alpha c {\Delta t}^2 \phi _1 (c (3 \alpha \ddot \phi _1
+2 c \phi _2^{(3)})
-2 \alpha^2 \dot \phi _2)
+(c^4 {\Delta t}^2
-{\Delta x}^2) (4 \phi _1^{(3)} \dot \phi _1
+4 \phi _2^{(3)} \dot \phi _2
+3 (\ddot \phi _1)^2
+3 (\ddot \phi _2)^2)\Big)
+\frac{h^4}{720} \Big(\alpha^6 {\Delta t}^4 \phi _1^2
+\alpha^6 {\Delta t}^4 \phi _2^2
+\alpha c {\Delta t}^4 \phi _1 (6 \alpha^4 \dot \phi _2
+c (
-15 \alpha^3 \ddot \phi _1
-20 \alpha^2 c \phi _2^{(3)}
+15 \alpha c^2 \phi _1^{(4)}
+6 c^3 \phi _2^{(5)}))
+\alpha c {\Delta t}^4 \phi _2 (
-6 \alpha^4 \dot \phi _1
-15 \alpha^3 c \ddot \phi _2
+20 \alpha^2 c^2 \phi _1^{(3)}
+15 \alpha c^3 \phi _2^{(4)}
-6 c^4 \phi _1^{(5)})
+(c^6 {\Delta t}^4
-{\Delta x}^4) (6 \phi _1^{(5)} \dot \phi _1
+6 \phi _2^{(5)} \dot \phi _2
+15 \phi _1^{(4)} \ddot \phi _1
+15 \phi _2^{(4)} \ddot \phi _2
+10 (\phi _1^{(3)})^2
+10 (\phi _2^{(3)})^2)\Big)
+\mathcal{O}(h^6)
		\end{autobreak},
	\end{align*}
	%
	where we have subtracted terms with a coefficient of the form $h^{2k+1}$ $(k\in \N_0)$ as these cannot have an impact on the dynamics because the discretisation scheme of the derivatives $u_{xx}$ and $u_{tt}$ is of even order. Indeed, they are in the kernel of the Euler operator $\mathcal {E}$ and do not influence the following computations. As before, $V$ is evaluated at $\|\phi\|^2 = \langle \phi,\phi \rangle$.
	
	We can explicitly verify that $\mathcal {E}{\mathcal L}$ recovers \eqref{eq:Taylor5PtStencilSymmetric} such that Palais' principle of symmetric criticality holds.
	As the highest derivative of $\phi$ in $\mathcal L$ is of order 5, we consider the coordinate $q = (q^1,q^2\ldots,q^5) = (\phi,\dot\phi,\ldots,\phi^{(4)})$ and compute the corresponding conjugate momenta $p = (p_1,\ldots,p_5)$ as
	\[
	p_i^j = \sum_{k=0}^{5-i} (-1)^k \frac{\d^{k}}{\d \xi^{k}}\frac{\p \mathcal{L}}{\p \phi_j^{(k+i)}}, \qquad j\in{1,2}.
	\]
	The Hamiltonian
	\[
	H = \sum_{i=1}^5 \langle \dot q^i , p_i \rangle - \mathcal L
	\]
	expressed in the dynamical variables $(\phi,\dot\phi,\ldots,\phi^{(4)},\phi^{(5)})$ is given as
	
	\begin{align*}
		\begin{autobreak}
			H=
\frac{1}{2} (-\alpha^2 \|\phi\|^2+(c^2-1)\|\dot \phi\|^2-V)
+\frac{h^2}{24}  \Bigg(\alpha^4 {\Delta t}^2\|\phi\|^2
-6 \alpha^2 c^2 {\Delta t}^2 \|\dot \phi\|^2
-8 \alpha c^3 {\Delta t}^2( \dot \phi _2 \ddot \phi _1-\dot \phi _1 \ddot \phi _2)
+2 (c^4 {\Delta t}^2-\Delta x^2) (\phi _1^{(3)} \dot \phi _1+\phi _2^{(3)} \dot \phi _2)
+(\Delta x^2-c^4 {\Delta t}^2 )\|\ddot \phi \|^2
\Bigg)
+\mathcal{O}(h^4)
		\end{autobreak},
	\end{align*}
	where the 4th order terms are computed but not presented here.
	
	\begin{remark}
		The Hamiltonian $H$ corresponds to the symplectic structure $\Omega= \sum_{i=1}^5 \sum_{j=1}^2\d q^i_j \wedge \d p_i^j$ where $q$ and $p$ are functions of $\{\phi\}_{5\mathrm{jet}}=(\phi,\dot\phi,\ldots,\phi^{(5)})$. For the computation of the modified system, it is not necessary to compute $\Omega$ explicitly. However, we remark that in the frame $\frac{\p}{\p \phi_1},\frac{\p}{\p \phi_2}\ldots,\frac{\p}{\p \phi_1^{(5)}},\frac{\p}{\p \phi_2^{(5)}}$ the skew symmetric matrix $\Omega^{\mathrm{MAT}}$ describing $\Omega$ is constant, i.e.\ does not depend on the base point $\{\phi\}_{5\mathrm{jet}}$. This is because $\mathcal L$ is a quadratic polynomial in $\{\phi\}_{5\mathrm{jet}}$. Hamilton's equations in $z=\{\phi\}_{5\mathrm{jet}}$ coordinates, i.e.\ $\dot z = (\Omega^{\mathrm{MAT}})^{-1}\nabla H(z)$, correspond to a first-order formulation of the truncated, high-order ODE \eqref{eq:Taylor5PtStencilSymmetric}.
	\end{remark}
	
	Now we express $H$ on-shell, i.e.\ we substitute second and higher order derivatives of $\phi$ with expressions in $(\phi,\dot\phi)$ obtained as in \eqref{eq:Phi2Elliminated}. We obtain the following expression:
	\begin{equation}\label{eq:Hmod}
		\begin{split}
			\mathcal H &=
			\frac{1}{2} ( -\alpha^2 \|\phi\|^2+(c^2-1)\|\dot \phi\|^2-V)\\
			&\phantom{=}+h^2
			\left(d_1 \|\phi\|^2
			+d_2 \|\dot \phi\|^2
			+ d_3 \langle \dot \phi, J \phi\rangle
			+ d_4 \langle \phi, \dot \phi \rangle^2
			\right)
			+\mathcal O(h^4)
		\end{split}
	\end{equation}
	
	\begin{align*}
		d_1 &=\frac{2 \alpha ^2 V' ({\Delta x}^2-c^4 {\Delta t}^2)+(V')^2 ({\Delta x}^2-c^4 {\Delta t}^2)+\alpha ^4 ((1-2 c^2) {\Delta t}^2+{\Delta x}^2)}{24 (c^2-1)^2}\\
		d_2&=\frac{\alpha ^2 (-3 c^4 {\Delta t}^2+c^2 (5 {\Delta x}^2-3 {\Delta t}^2)+{\Delta x}^2)+(c^2-1) V' (c^4 {\Delta t}^2-{\Delta x}^2)}{12 (c^2-1)^2}\\
		d_3&=\frac{\alpha  c (c^2 {\Delta t}^2-{\Delta x}^2) (\alpha ^2+V')}{3 (c^2-1)^2}\\
		d_4&=\frac{V'' (c^4 {\Delta t}^2-{\Delta x}^2)}{6 (c^2-1)}
	\end{align*}
	The above Hamiltonian has been expressed in invariant terms of the rotation symmetry described in \cref{rem:RotInvExact,rem:symDiscr}.
	Higher order terms are reported in \cite{multisymplecticSoftware}.

	The symplectic structure is given as $\omega = \sum_{i=1}^5 \sum_{j=1}^2\d q^i_j \wedge \d p_i^j$ where $q$ and $p$ are functions of $(\phi,\dot\phi)$. In the frame $\frac{\p}{\p \phi_1},\frac{\p}{\p \phi_2},\frac{\p}{\p \dot\phi_1},\frac{\p}{\p \dot \phi_2}$ the 2-form $\omega$ is represented by the following matrix.\footnote{In  Mathematica the matrix $\omega^{\mathrm{MAT}}$ can be computed by applying the operator {\tt TensorWedge} to the gradients of $q^i_j(\phi,\dot\phi)$ and $p^i_j(\phi,\dot\phi)$ and summing over the indices.}
	
	\[\omega^{\mathrm{MAT}}=
	\begin{pmatrix}
		0 & 2 \alpha c & 1-c^2 & 0 \\
		-2 \alpha c & 0 & 0 & 1-c^2 \\
		c^2-1 & 0 & 0 & 0 \\
		0 & c^2-1 & 0 & 0 \\
	\end{pmatrix}
	+h^2
	\begingroup 
	\renewcommand{\arraystretch}{2}
	\begin{pmatrix}
		w_1J&Z\\-Z&w_2J
	\end{pmatrix}
	\endgroup
	+\mathcal O(h^4)
	\]
	
	\begin{align*}
		w_1& = \frac{\alpha c}{3 (c^2-1)}
		\big( \alpha(\Delta x^2-\Delta t^2) + (\Delta x^2-2c^2\Delta t^2+c^4 \Delta t^2)(V'+\|\phi\|^2 V'')\big)\\
		w_2&=\frac{\alpha c(c^2 \Delta t^2-\Delta x^2)}{3 (c^2-1)}\\
		Z &=\left( -\frac{\alpha^2 (c^2 ((c^2-3) {\Delta t}^2+{\Delta x}^2)+{\Delta x}^2)}{6 (c^2-1)^2}+\frac{(c^2-1)(\Delta x^2-c^4\Delta t^2)}{6 (c^2-1)^2}V'\right)
		\begin{pmatrix}
			0&1\\1&0
		\end{pmatrix}\\
		&\phantom{=}-\frac{(c^4 \Delta t^2-\Delta x^2)}{3 (c^2-1)}V''
		\begin{pmatrix}
			\phi_1^2&\phi_1\phi_2\\
			\phi_1\phi_2&\phi_2^2
		\end{pmatrix}\\
		J&= \begin{pmatrix}
			0&1\\-1&0
		\end{pmatrix}
	\end{align*}
	The fourth-order expressions are quite long and are omitted here. We remark, however, that while the second-order terms do not explicitly depend on $\dot \phi$, the fourth-order terms do.
	Hamilton's equations for $(\mathrm{Jet}^1(Q),\omega,\mathcal H)$ recover\footnote{When computing in local coordinates $\dot z = (\omega^{\mathrm{MAT}})^{-1}\nabla H(z)$ with $z=(\phi,\dot \phi)$ the matrix inversion $(\omega^{\mathrm{MAT}})^{-1}$ needs to be computed in the sense of formal power series.} \eqref{eq:Phi2Elliminated} with an $\mathcal{O}(h^6)$ error.
	
	\begin{remark}\label{rem:NoLStandard}
		In the general case, the symplectic structure $\omega$ contains $\d \dot \phi_i \wedge \d \dot\phi_j$-terms, i.e.\ the fibres of $\mathrm{Jet}_1(Q)$ spanned by the vertical vector fields $\frac{\p}{\p \dot\phi_1},\frac{\p}{\p \dot\phi_2}$ are {\em not} Lagrangian with respect to $\omega$. Therefore, there is no primitive $\lambda$ with $\d \lambda=\omega$ which does {\em not} involve any $\d \dot \phi_j$-terms, i.e.\ vanishes on vertical vector fields.
		Thus, the Lagrangian density $\lambda - \mathcal H \d \xi$ does {\em not} admit the form $L(\phi,\dot \phi)\d \xi$ for any $\lambda$ and $L$ because $L$ will depend on $\ddot \phi$.
		This means, if one wants to express the system as a first-order variational system in the familiar form of first-order variational principles \eqref{eq:AssFormL}, then a change of coordinates $(\phi,\dot \phi)\mapsto( q, \dot {  q})$ is necessary, such that the distribution $\mathcal D$ spanned by $\frac{\p}{\p \dot {  q}^1},\frac{\p}{\p \dot {  q}^2}$ is Lagrangian for $\omega$. The existence of a 1-form $\lambda$ with kernel $\mathcal D$ is then guaranteed by \cite[Cor.\ 15.7]{Libermann1987}.
		By the existence of Darboux coordinates, coordinates can be chosen such that $\lambda =  \dot { q}^1 \d { q}^1+\dot { q}^2 \d { q}^2$. However, one cannot expect the coordinate transformation to admit a closed form in the general case.
		(See also \cref{rem:NoBundle}.)
	\end{remark}

	\begin{remark}
		The Lagrangian $\mathcal{L}$ has been obtained from $L_\Delta$ by a series expansion. Indeed, for the computations $\mathcal{L}$ could be substituted by any Lagrangian which induces equivalent Euler-Lagrange equations. This is in contrast to the backward error analysis methods for classical variational integrators developed in \cite{Vermeeren2017}, where the Euler-McLaurin formula is used to translate the discrete action to an integral expression over a Lagrangian given as a formal power series.
	\end{remark}
	
	Following \cref{prop:presConsQ} we can compute the modified conserved quantity resulting from the rotational invariance of the original equation and the five-point stencil. If $M$ is the order of the highest derivative in the truncated Lagrangian $\mathcal{L}$ for some truncation index $N$ then by Noether's theorem the conserved quantity is given as
	\begin{align*}
		I_{\mathrm{rot}}^{\Delta}(\phi,\dot \phi,\ldots,\phi^{(M)})
		= \sum_{m=1}^M \sum_{k=0}^{m-1}(-1)^k \left\langle \frac{\d^k}{\d \xi^k} \nabla_{\phi^{(m)}}{\mathcal L}, J^T \phi^{(m-1-k)}\right\rangle.
	\end{align*}
	This yields the reduced quantity
	\begin{align*}
		I_{\mathrm{rot}}^{\mathrm{mod}}(\phi,\dot \phi)
		&= \alpha c \|\phi\|^2 + (c^2-1) \langle J \dot \phi, \phi \rangle\\
		&\phantom{=}+\frac{h^2}{6 \left(c^2-1\right)^2}\left(b_1\|\phi\|^2+b_2\|\dot \phi\|^2+b_3 \langle J\dot \phi,\phi\rangle\right)+\mathcal O(h^4)
	\end{align*}
	with
	\begin{align*}
		b_1&=\alpha  c \left(\alpha ^2 \left({\Delta x}^2-{\Delta t}^2\right)+V' \left(c^2 \left(c^2-2\right) {\Delta t}^2+{\Delta x}^2\right)\right)\\
		b_2&= \alpha  c \left(c^2-1\right) \left(c^2 {\Delta t}^2-{\Delta x}^2\right)\\
		b_3 &=\alpha ^2 \left(c^4 {\Delta t}^2+c^2 \left({\Delta x}^2-3 {\Delta t}^2\right)+{\Delta x}^2\right)+\left(c^2-1\right) V' \left(c^4 {\Delta t}^2-{\Delta x}^2\right).
	\end{align*}

	\begin{figure}
		\centering
		\begin{subfigure}{0.45\textwidth}
			\centering
			\includegraphics[width=\textwidth]{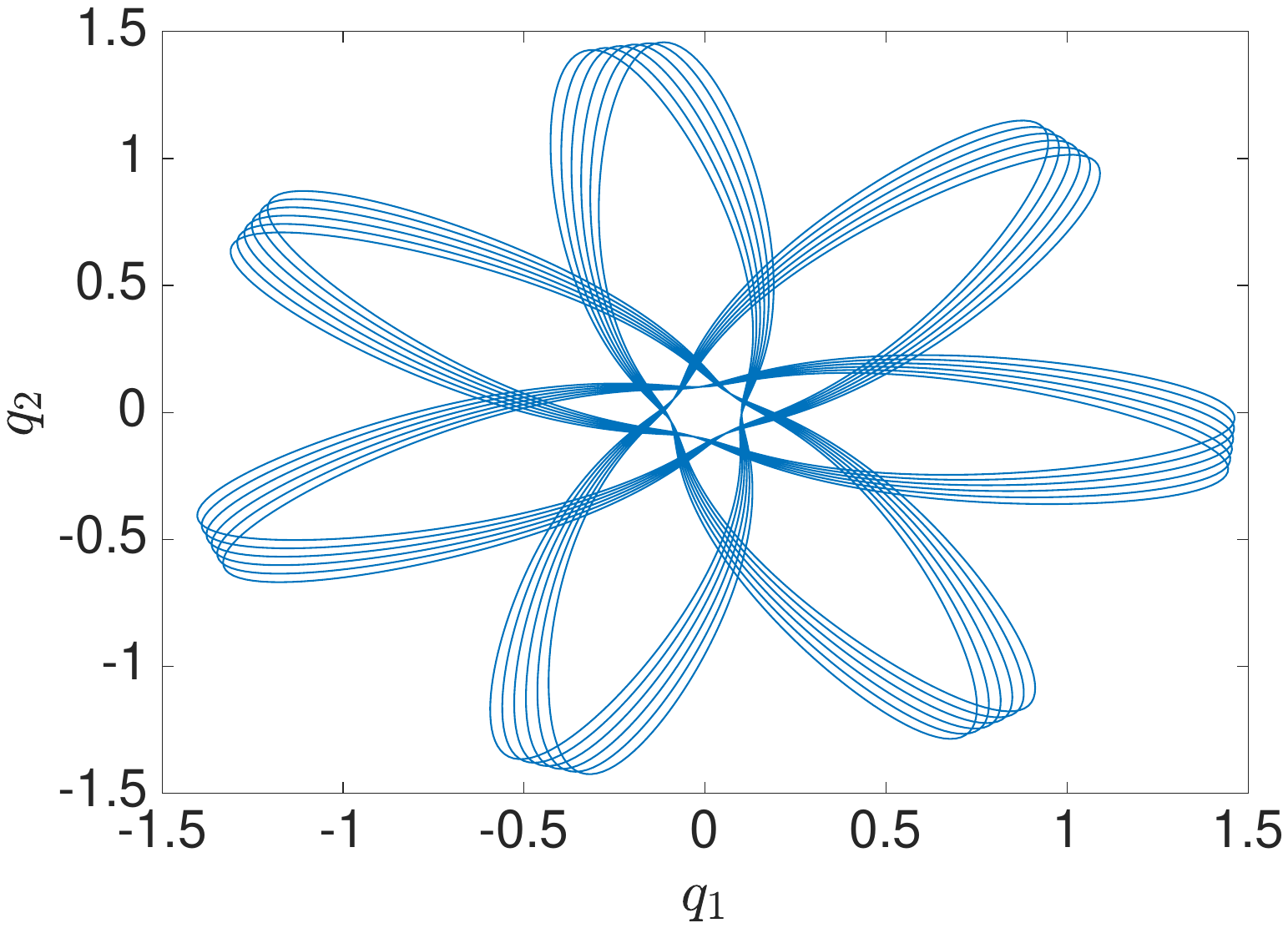}
			\subcaption{Phase plot, $\xi\in[0, 500]$}
			
		\end{subfigure}
		\begin{subfigure}{0.49\textwidth}
			\centering
			\includegraphics[width=\textwidth]{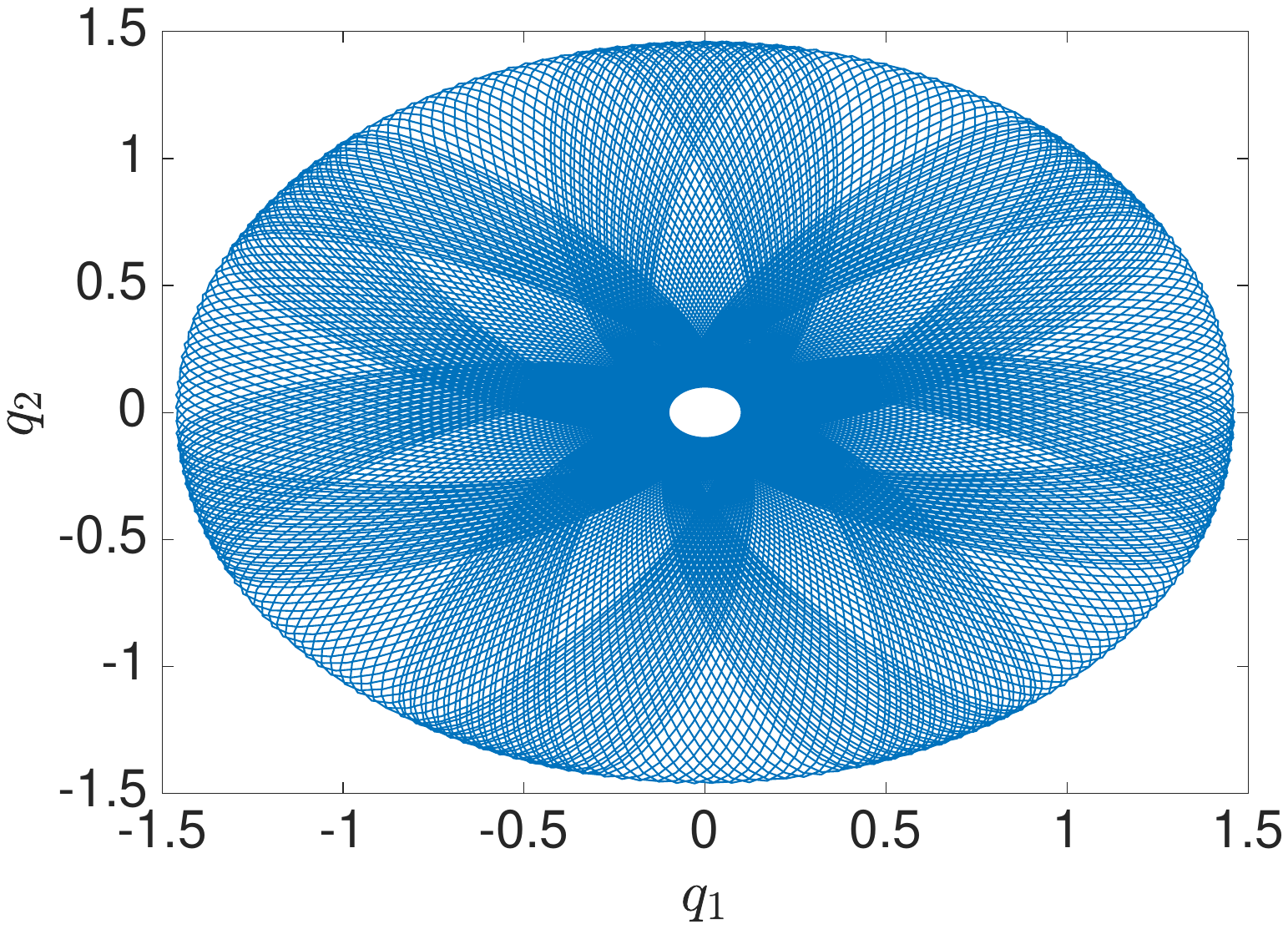}
			\subcaption{Phase plot, $\xi\in[0,4000]$}
			
		\end{subfigure}
		\begin{subfigure}{0.49\textwidth}
			\centering
			\includegraphics[width=\textwidth]{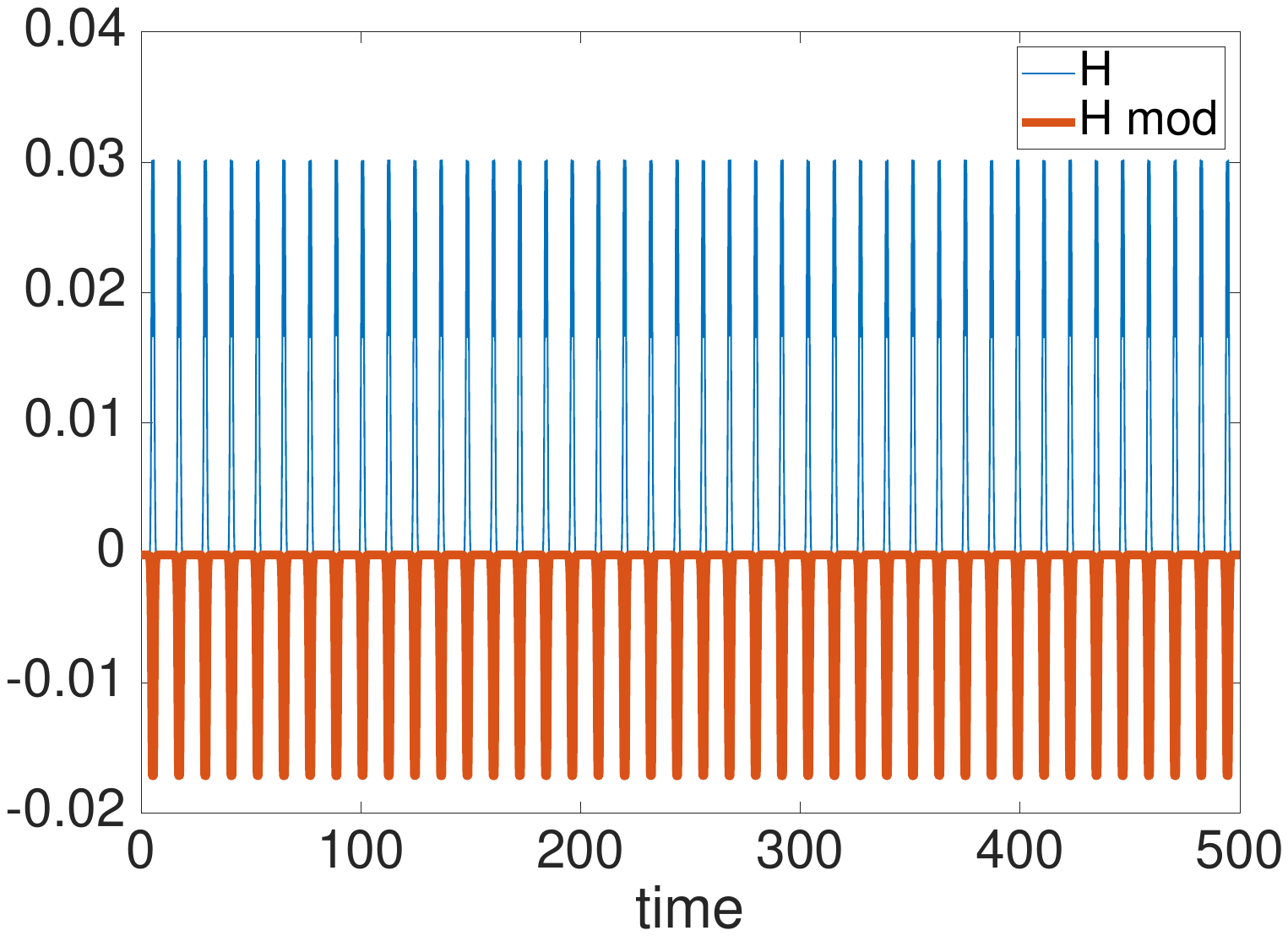}
			\subcaption{$H$ (blue) and $\mathcal H$ (red)}
			
		\end{subfigure}
		\begin{subfigure}{0.49\textwidth}
			\centering
			\includegraphics[width=\textwidth]{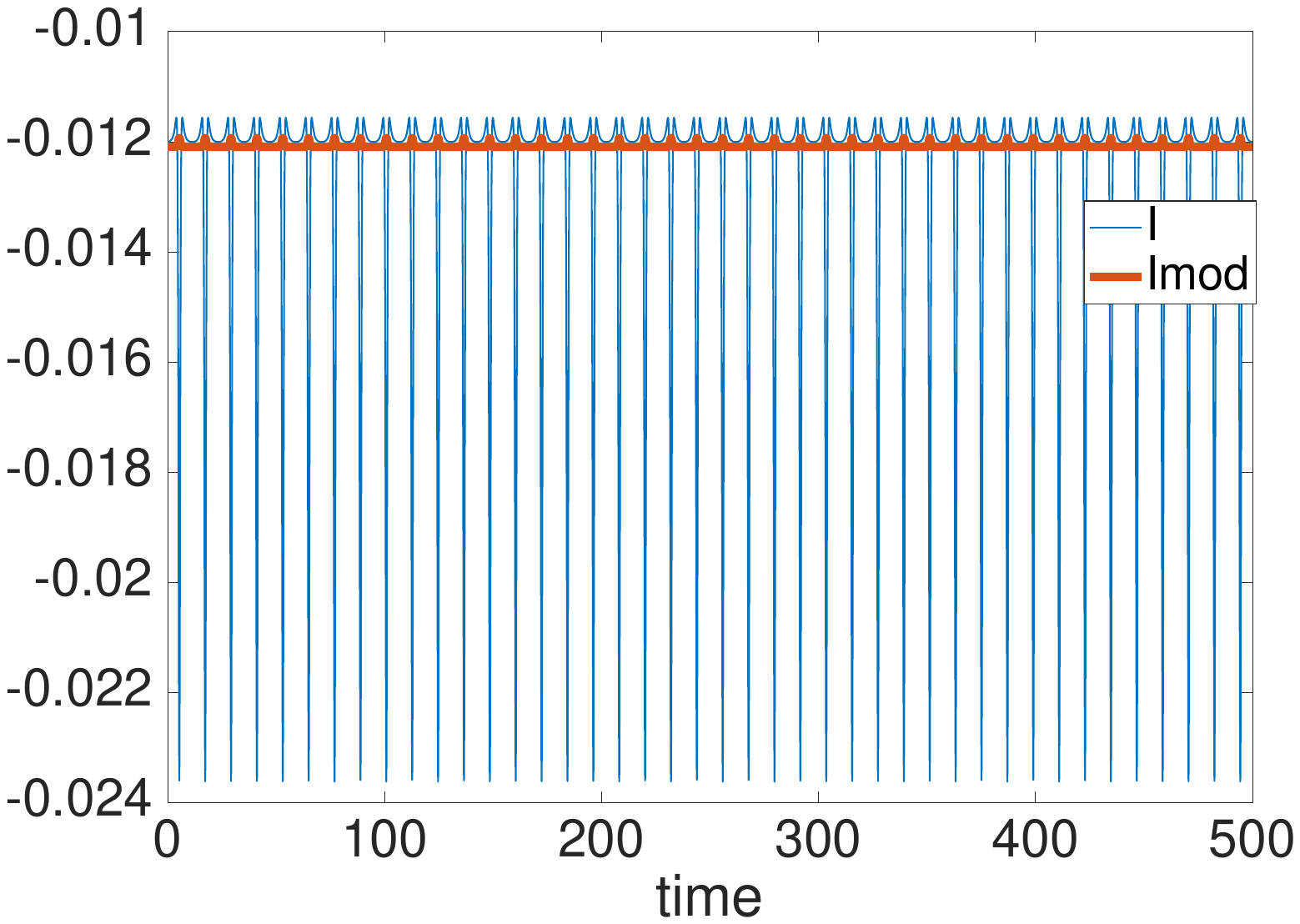}
			\subcaption{$I_{\mathrm{rot}}$ (blue) and $I_{\mathrm{rot}}^{\mathrm{mod}}$ (red)}
			
		\end{subfigure}
		\caption[Dynamics modified equation]{Numerical integration of the ODE \eqref{eq:Phi2Elliminated} truncated after $\mathcal O(h^2)$ terms with $V(s) = -0.1 s^4 +s$, $\Delta x= 0.1$, $\Delta t = 0.15$, $\alpha = 0.3$, $c = 2$.
			All numerical computations have been performed in the Darboux variables $(\mathfrak q, \mathfrak p)$ of the continuous system using the implicit midpoint rule combined with fixed-point iterations. Therefore, the integration is symplectic modulo second order terms. The plots show a phase plot of a motion initialised at $(\mathfrak q, \mathfrak p)=(-0.11,-0.01,-0.1,0.1)$ and the behaviour of the Hamiltonian $H$ of the exact system and the Hamiltonian $\mathcal H$ of the modified system truncated after $\mathcal O(h^2)$ terms as well as the behaviour of the conserved quantity of the exact system $I_{\mathrm{rot}}$ and of the modified system $I_{\mathrm{rot}}^{\mathrm{mod}}$ truncated after $\mathcal O(h^2)$ terms along the motion.
		}\label{fig:ExpModSystem}
	\end{figure}
	
		\begin{figure}
		\centering
		\begin{subfigure}{\textwidth}
			\centering
			\includegraphics[width=0.8\textwidth]{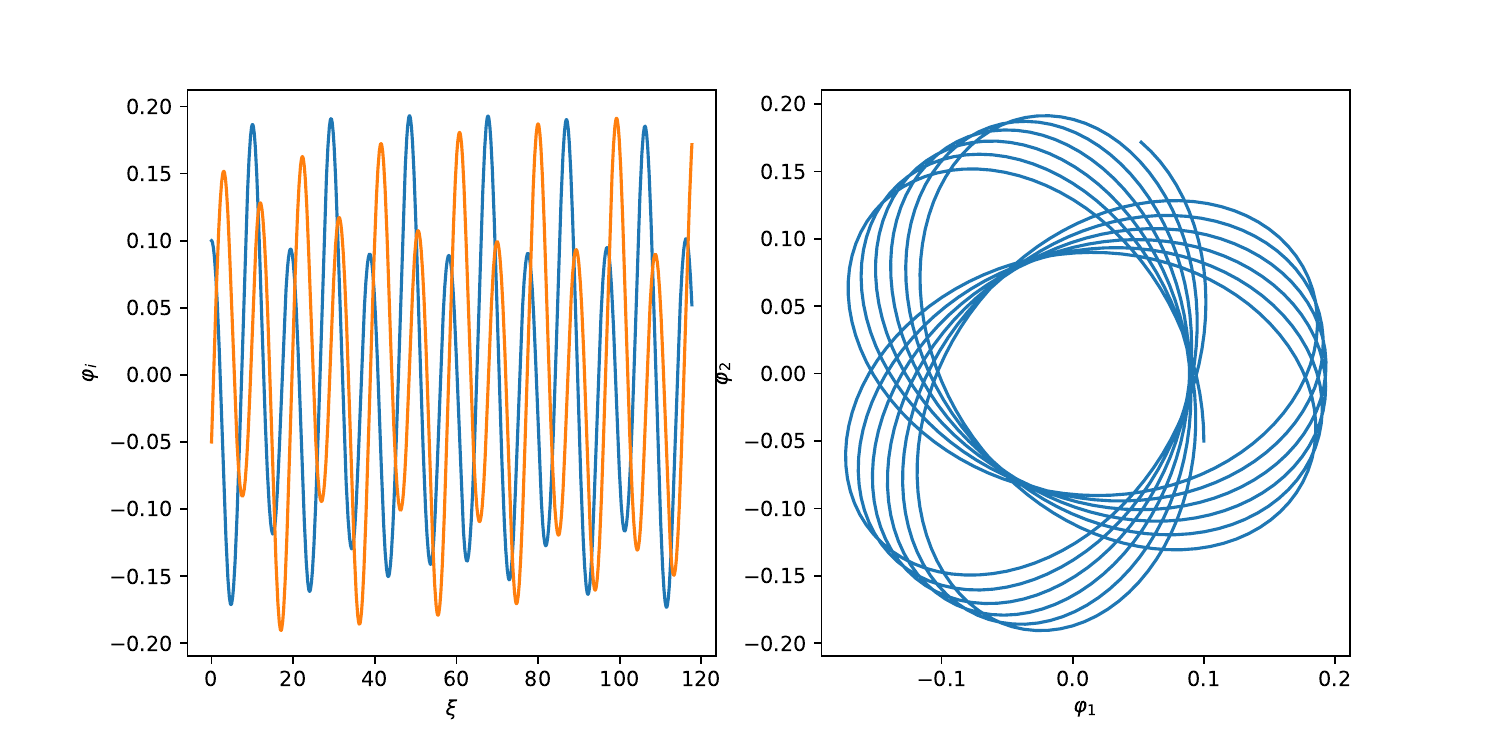}
			\subcaption{dynamics and phase plot, $\xi\in[0,120]$}
			
		\end{subfigure}
		\begin{subfigure}{0.49\textwidth}
			\centering
			\includegraphics[width=\textwidth]{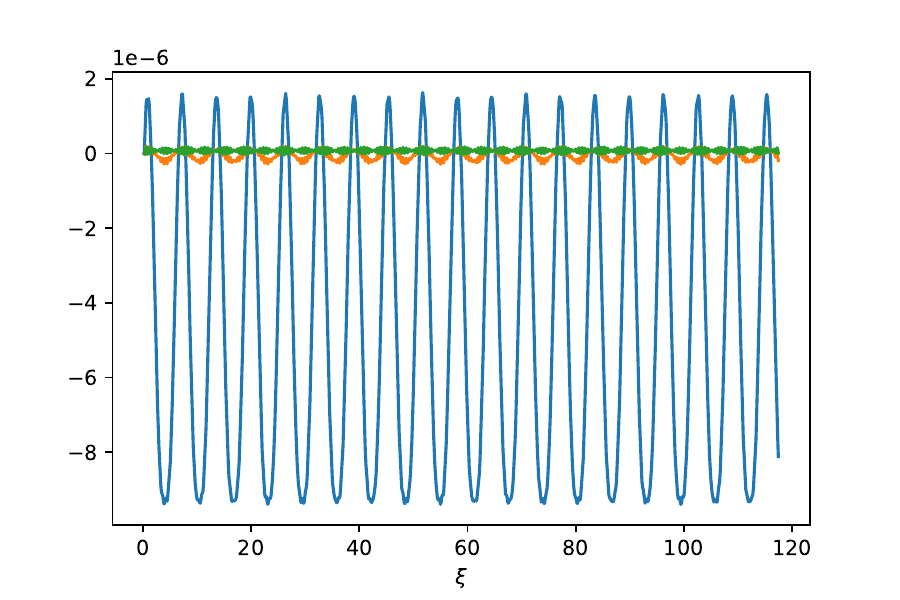}
			\subcaption{$\mathcal{H}$ truncated to 0th (blue), 2nd (orange), 4th (green) order}
			
		\end{subfigure}
		\begin{subfigure}{0.49\textwidth}
			\centering
			\includegraphics[width=\textwidth]{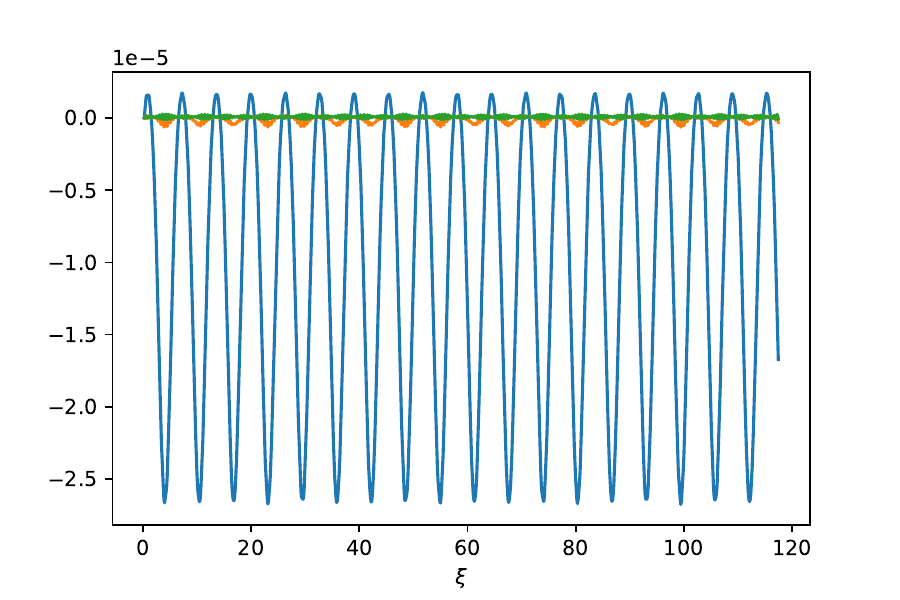}
			\subcaption{$I_{\mathrm{rot}}^{\mathrm{mod}}$ to 0th (blue), 2nd (orange), 4th (green) order}
			
		\end{subfigure}
		\caption[Dynamics functional equation]{
			When $c \Delta t/\Delta x$ is rational, the functional equation \eqref{eq:5PtStencilSymmetric} yields a multistep formula. The series parameter $h$ is set to 1. We use $V(s)=s^2$, $\Delta t = 0.15$, $\Delta x = 2 c \Delta t$. Let $\Delta \tau = c \Delta t$. To initialise the scheme, values at $\xi = \Delta \tau, 2\Delta \tau, 3\Delta \tau$ are obtained by integrating \eqref{eq:Phi2Elliminated} truncated to 4th order with high accuracy with the initial condition $(\phi(0),\dot{\phi}(0)) = ((0.1,-0.05),(0,0.1))$.
		}\label{fig:Stepping}
	\end{figure}

	The modified system $(\mathrm{Jet}^1(Q),\omega,\mathcal H)$ is completely Liouville integrable up to any order in the series parameter $h$.
	The plot of a motion of the second-order truncated system and the behaviour of $H$ and $\mathcal H$ for example data and the behaviour of the quantities $I_{\mathrm{rot}}$ and $I_{\mathrm{rot}}^{\mathrm{mod}}$ truncated after $\mathcal O(h^2)$ terms can be seen in \cref{fig:ExpModSystem}. The plotted system is $\mathcal O(h^4)$ close to a completely integrable system.
	Similarly, \cref{fig:Stepping} shows analogous experiments but the motion is obtained by using \eqref{eq:5PtStencilSymmetric} as a multistep formula
	when the ratio $c \Delta t/\Delta x$ is rational. For consistent initialisation of the multistep formula the modified equation \eqref{eq:Phi2Elliminated} was used truncated to 4th order to eliminate effects of spurious solutions.
	The effects of spurious solutions can be amplified if $V$ contains terms of high polynomial order, the rational relation of $\Delta x$ and $c \Delta t$ requires large denominators such that the multistep formula needs to be initialised at many steps, or a low order truncation of the modified equation is used in the initialisation process.


	\subsection{Special cases}\label{sec:cases}
	
	Following up on \cref{rem:NoLStandard}, we consider the following special cases for which the symplectic structure $\omega$ does not contain $\d \dot \phi_i \wedge \d \dot\phi_j$-terms, i.e.\ the fibres of $\mathrm{Jet}_1(Q)$ spanned by the vertical vector fields $\frac{\p}{\p \dot\phi_1},\frac{\p}{\p \dot\phi_2}$ are Lagrangian with respect to $\omega$.
	
	\begin{itemize}
		\item
		$\alpha=0$ (non-rotating travelling wave),
		\item
		$c=0$ (standing wave),
		\item
		$\Delta x = c \Delta t$ (step sizes fulfil a special relation)
	\end{itemize}
	
	This can be seen up to the computed order by inspecting the bottom right $2 \times 2$ submatrix in $\omega^{\mathrm{MAT}}$.
	By \cite[Cor.\ 15.7]{Libermann1987} there exists a 1-form $\lambda$ with $\d \lambda = \omega$ such that the vertical vector fields span its kernel. In other words, the 1-form $\lambda$ does not contain any $\d \dot\phi_1$ or $\d \dot\phi_2$ terms. Thus, there exists a first-order Lagrangian in the original variables $L(\phi,\dot \phi)$ such that the Euler-Lagrange equations recover the dynamics.
	Indeed, this can be shown independently of our computation up to any order for the cases $c=0$ and $\Delta x = c \Delta t$ because in these cases the functional equation \eqref{eq:5PtStencilSymmetric} relates three equally-spaced points rather than five unequally spaced ones. This is argued in the following.
	(The case $\alpha=0$ is treated in Section \ref{sec:MSApproach}.)
	
	For $c=0$ and $h=1$ the functional equation \eqref{eq:5PtStencilSymmetric} reduces to
	\begin{equation}\label{eq:c0ODE}
		-\frac{\phi(\xi+\Delta x)-2\phi(\xi)+\phi(\xi-\Delta x)}{\Delta x^2}
		-\nabla W(\phi(\xi))
		+\frac{R(-\Delta t)-2 I + R(\Delta t)}{\Delta t^2} \phi(\xi)=0,
	\end{equation}
	where $W(\phi) = \frac 12 V(\|\phi\|^2)$.
	The relation \eqref{eq:c0ODE} arises as the discrete Euler-Lagrange equations
	\begin{equation}\label{eq:discrEL}
		D_1L_{\Delta}(\phi(\xi),\phi(\xi+\Delta x))+D_2L_{\Delta}(\phi(\xi-\Delta x),\phi(\xi))=0
	\end{equation}
	for
	\[
	L^{c=0}_{\Delta}(a,b) = \frac {\|a-b\|^2}{2 \Delta x^2}-W(a)+\frac 1{2\Delta t^2}a^T (R(-\Delta t)-2I + R(\Delta t)) a,
	\]
	where $D_1$ denotes the derivative with respect to the variables passed into the first argument of $L_\Delta$ and $D_2$ denotes the derivative with respect to the second argument.
	The functional equation \eqref{eq:c0ODE} is a discretisation of
	\[
	-\ddot \phi - \nabla W(\phi)+\frac 1 {\Delta t^2} (R(-\Delta t)-2 I + R(\Delta t))\phi=0
	\]
	which is the Euler-Lagrange equation to the Lagrangian
	\[
	L^0_{c=0} =\frac 12 \|\dot \phi\|^2 -W(\phi)+\frac 1 {2 \Delta t^2} \phi^\top (R(-\Delta t)-2 I + R(\Delta t))\phi.
	\]
	By \cite{Vermeeren2017} there exists a modified Lagrangian given as the formal power series
	\[
	L_{c=0} = L_{c=0}^0(\phi,\dot{\phi}) + h^2 L_{c=0}^2(\phi,\dot{\phi}) + h^4 L_{c=0}^2(\phi,\dot{\phi})+\ldots.
	\]
	The expression can either be found up to any order in $h$ using the techniques from \cite{Vermeeren2017} or with an educated guess for $L_{c=0}$, see \begin{small}{\tt Computational\_Results\_documented.pdf}\end{small} in \cite{multisymplecticSoftware}.

	The situation is similar for $\Delta x = c \Delta t$. The functional equation \eqref{eq:5PtStencilSymmetric} reduces to
	\begin{equation}\label{eq:dxcdt}
		\begin{split}
			\frac{c^2}{\Delta x^2}(R(-\Delta t)\phi(\xi +\Delta x) - 2\phi(\xi) + R(\Delta t) \phi(\xi-\Delta x) )\\
			- \frac {\phi(\xi + \Delta x)-2\phi(\xi)+\phi(\xi - \Delta x)}{\Delta x^2}-\nabla W(\phi)=0.
		\end{split}
	\end{equation}
	The functional equation \eqref{eq:dxcdt} arises as the discrete Euler--Lagrange equations \eqref{eq:discrEL} for the discrete Lagrangian
	\[
	L_\Delta^{\Delta x = c \Delta t}(a,b) = -\frac {c^2}{2\Delta x^2} \left\|R\left(-\frac 12 \Delta t\right)b-R\left(\frac 12 \Delta t\right)a\right\|^2 + \frac 1 {2 \Delta x^2}\|b-a\|^2-W(a).
	\]
	Equation \eqref{eq:dxcdt} is a discretisation of
	\[
	(c^2-1)\ddot \phi -2 \alpha c J \dot \phi - \alpha^2\phi-\nabla W(\phi) =0
	\]
	which is governed by
	\begin{equation}\label{eq:L0dxcdt}
		L_{\Delta x = c \Delta t}^0(\phi,\dot \phi) =\frac{1}{2} \left( \alpha^2 \langle \phi,\phi\rangle - 2 \alpha c \langle J \phi, \dot{\phi}\rangle + (c^2-1) \langle \dot{\phi},\dot{\phi}\rangle\right) + W(\phi),
	\end{equation}
	i.e.\ we recover $L^0$ from \eqref{eq:L0Def}.
	In analogy to the case $c=0$, up to any order in $h$ there exists a modified first-order Lagrangian in the original variables $\phi,\dot \phi$ such that its Euler--Lagrange equations recover \eqref{eq:Phi2Elliminated}.
	Computational results are reported in \cite{multisymplecticSoftware}.



	\section{Computation of modified Lagrangian via P-series ansatz}\label{sec:MSApproach}
	
	We now turn to the case of non-rotating travelling waves for the nonlinear wave equation, for
	which the discrete travelling wave equation (see \eqref{eq:dtwe}) is
	\begin{equation*}
		\begin{aligned}
			0&=\frac{1}{{\Delta t}^2} \left( \phi(\xi+c\Delta t) - 2 \phi(\xi) + \phi(\xi-c\Delta t) \right) \\
			&- \frac{1}{{\Delta x}^2} \left(  \phi(\xi+\Delta x) - 2 \phi(\xi) + \phi(\xi-\Delta x) \right) \\
			&-\nabla W(\phi(\xi))
		\end{aligned}
	\end{equation*}
	where $\phi\colon\R\to Q=\R^d$ and $W\colon \R^d\to\R$. We know from \cref{thm:ConsPalais} that $\phi$ has a modified equation to all orders with Hamiltonian
	structure $(\mathrm{Jet}_1(Q),\omega,\mathcal{H})$.
		\begin{figure}
		\begin{center}
			\begin{tikzpicture}[ scale=0.85]
				
				\draw[thick, -] (0.5,0) -- (15.5,0);
				
				\draw (1 cm, 8pt) -- (1 cm, -8pt) node[anchor = north] {$\hat \xi$};

				\draw (2 cm, 3pt) -- (2 cm, -3pt);
				\draw (3 cm, 3pt) -- (3 cm, -3pt);
				\draw (4 cm, 8pt) -- (4 cm, -8pt) node[anchor = north] {$\hat \xi + \left(\frac 12 - \frac m{2n}\right)\Delta s$};
				\draw (5 cm, 3pt) -- (5 cm, -3pt);
				\draw (6 cm, 3pt) -- (6 cm, -3pt);
				\draw (7 cm, 3pt) -- (7 cm, -3pt);
				
				\draw (8 cm, 8pt) -- (8 cm, -8pt) node[anchor = north] {$\hat \xi + \frac 12 \Delta s$};
				\draw (9 cm, 3pt) -- (9 cm, -3pt);
				\draw (10 cm, 3pt) -- (10 cm, -3pt);
				\draw (11 cm, 3pt) -- (11 cm, -3pt);
				\draw (12 cm, 8pt) -- (12 cm, -8pt) node[anchor = north] {$\hat \xi + \left(\frac 12+\frac{m}{2n}\right)\Delta s$};
				\draw (13 cm, 3pt) -- (13 cm, -3pt);
				\draw (14 cm, 3pt) -- (14 cm, -3pt);
				\draw (15 cm, 8pt) -- (15 cm, -8pt) node[anchor = north] {$\hat \xi + \Delta s$};
				
				\draw (9 cm, 8pt) -- (9 cm, 4pt) ;
				\draw (10 cm, 8pt) -- (10 cm, 4pt) ;
				\draw (9 cm, 6pt) -- (10 cm, 6pt) ;
				\draw (9.5 cm,6pt) node[anchor = south] {$\frac{\Delta s}{2n}$};
				
			\end{tikzpicture}

		\end{center}
		\caption{
			Interpretation of \eqref{eq:dtwe} as a multistep formula for $\frac mn = \frac{\Delta x}{c\Delta t}<1$.
			The variable $\hat \xi$ corresponds to $\xi - c \Delta t$ when comparing with \eqref{eq:dtwe} and $\Delta s = 2 c \Delta t$.
		}\label{fig:MSFormulaIlluIrr}
	\end{figure}
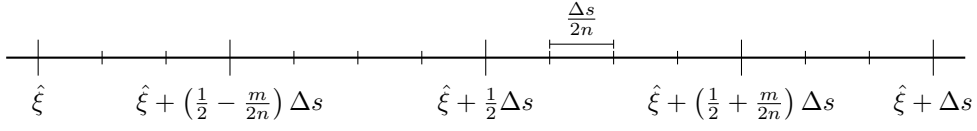
	However, we can also interpret the above equation as a linear multistep formula, see \cref{fig:MSFormulaIlluIrr}. By our discussion of blended backward error analysis for linear multistep schemes \cite{BEAMulti}, the modified equation is equivalent the Euler-Lagrange equations for a 2nd order modified Lagrangian given as an $S$-series applied to a $P$-series in the original variables $(\phi,\dot{\phi})$.
	The structural information helps to quickly compute terms of the modified Lagrangian using an educated guess (ansatz).
	For an introduction to $P$-series, which arise in partitioned Runge--Kutta methods and multistep
	methods, see, for instance, \cite[III.2.1]{GeomIntegration}.

	\begin{proposition}\label{prop:Phi2SeriesVariational}
		The discrete travelling wave equation \eqref{eq:dtwe} has a second-order modified equation
		that arises as the Euler--Lagrange equation of a formal series
		\[
		L(\phi,\dot\phi) = \sum_{j=0}^\infty h^{2j}L^{2j}(\phi,\dot\phi).
		\]
		Moreover,
		\[L^0 = \frac{1}{2}  (c^2-1) \langle \dot{\phi},\dot{\phi}\rangle + W( \phi )
		\]
		and each term $L^{2j}$ for $j>0$ is a sum of elementary differentials of $W$ and $\dot\phi$ associated
		to bicoloured trees in which the leaves are black or white, the other nodes are black, and the
		sum of the degrees of the black nodes is $2j$.
	\end{proposition}
	
	\begin{proof}
		The statement follows from the authors' discussion of backward error analysis techniques for multistep methods \cite{BEAMulti}, in particular Remark 6.3, which is based on \cite{Chartier2006}: indeed, the multistep framework applies even for irrational quotients $\frac{c \Delta t}{\Delta x}$ (in which case \eqref{eq:dtwe} does not correspond to a practical multistep formula) because we can
		considering the non-polynomial characteristic functions
		\begin{align*}
			\rho(\tau) &=
			4 c^2
			-4\frac{c^2 \Delta t^2}{\Delta x^2}\tau^{\frac 12 - \frac{\Delta x}{2c \Delta t}}
			+8\left(\frac{c^2 \Delta t^2}{\Delta x^2} - c^2\right)\tau^{\frac 12}\\
			&-4\frac{c^2 \Delta t^2}{\Delta x^2}\tau^{\frac 12 + \frac{\Delta x}{2c \Delta t}}
			+4c^2 \tau\\
			\sigma(\tau) &= \tau^{\frac 12}
		\end{align*}
		which successfully recover the stencil \eqref{eq:dtwe} as \[\rho(e^D)\phi = \Delta s^2\sigma(e^D) \nabla W(\phi).\] Here the stepsize ${\Delta s}$ is given as ${\Delta s}= 2c\Delta t$ and $c\Delta t > \Delta x$. The case $c\Delta t < \Delta x$ is argued analogously.
		%
		%
		%
		%
	\end{proof}
	

	Elementary differentials and trees are motivated by their
	appearance in the modified Hamiltonians of partitioned Runge--Kutta methods \cite[III.2.1]{GeomIntegration}.
	The second order trees are
	\scalebox{0.35}{
		\begin{forest}
			for tree={fill=black}
			[[]]
	\end{forest}},
	with elementary differential $F_{2,1} = \sum_{i=1}^d W_i W_i$ and
	\scalebox{0.35}{
		\begin{forest}
			[[][]]
	\end{forest}}, with elementary differential
	$F_{2,2} = \sum_{i,j=1}^d W_{i,j}\dot\phi_i\dot\phi_j$.
	Here and in the following, subscripts indicate partial derivatives, e.g. $W_{ij} = \frac{\partial^2 W}{\partial \phi_i\partial\phi_j}$.
	There are four order 4 trees given as
	\begin{center}
		\scalebox{0.44}{
			\begin{forest}
				[[][][][]]
			\end{forest}\quad
			\begin{forest}
				[
				[, fill=black
				[][]
				]
				]
			\end{forest}\quad
			\begin{forest}
				for tree={fill,for descendants={fill=white}}
				[
				[, fill=black
				[]
				]
				[]
				]
			\end{forest}\quad
			\begin{forest}
				for tree={fill=black}
				[[[]]]
		\end{forest}}
	\end{center}
	with elementary differentials
	\[
	F_{4,1}=\sum_{i,j,k,l=1}^d W_{ijkl} \dot \phi_i \dot \phi_j \dot \phi_k\dot \phi_l, \quad
	F_{4,2}=\sum_{i,j,k=1}^d W_i W_{ijk}  \dot \phi_j \dot \phi_k,
	\]
	\[
	F_{4,3}=\sum_{i,j,k=1}^d W_{ij} \dot \phi_i W_{jk} \dot \phi_k, \quad
	F_{4,4}=\sum_{i,j=1}^d W_i W_{ij} W_j,
	\]
	respectively.
	There are 10 order 6 trees which we list below.\label{order6treesCompExMultiSympl}
	
	\scalebox{0.44}{
		\begin{forest}
			[[][][][][][]]
		\end{forest}\quad
		\begin{forest}
			[[,fill=black[][][]]]
		\end{forest}\quad
		\begin{forest}
			[[,fill=black[,fill=black][]][]]
		\end{forest}\quad
		\begin{forest}
			[[,fill=black[,fill=black[]]][]]
		\end{forest}
		
		\begin{forest}
			for tree={fill=black}
			[[[[]]]]
		\end{forest}\quad
		\begin{forest}
			[[,fill=black[,fill=black[][]]]]
		\end{forest}\quad
		\begin{forest}
			[[,fill=black[][][,fill=black]]]
		\end{forest}
		
		\begin{forest}
			[[,fill=black[]][][][]]
		\end{forest}\quad
		\begin{forest}
			[[,fill=black[][]][][]]
		\end{forest}\quad
		\begin{forest}
			for tree={fill=black}
			[[[][]]]
		\end{forest}
	}
	
	Thus \Cref{prop:Phi2SeriesVariational} states that the terms in the modified Lagrangian take the form $L^{2j} = \sum_k a_{j,k}F_{j,k}(\phi,\dot\phi)$ for suitable
	coefficients $a_{j,k}$.

	We compute the Euler-Lagrange equations for the ansatz described in \cref{prop:Phi2SeriesVariational} and solve for $\ddot{\phi}(\xi)$ using a series ansatz. Comparing coefficients with \eqref{eq:Phi2Elliminated} yields
	%
	%
	%
	\begin{align*}
		a_{2,2}&= 2 a_{2,1} = \frac{1}{12}(c^2-1)^{-1}\big(c^4 {\Delta t}^2-{\Delta x}^2\big),\\
		a_{4,1}&=(c^2-1)^{-3}b_1,\
		a_{4,2}=6(c^2-1)^{-4}b_1,\
		a_{4,3}= (c^2-1)^{-3}b_2,\\
		a_{4,4}&=3(c^2-1)^{-5}b_1,\
		a_{6,1} = (c^2-1)^{-3} b_3,\
		a_{6,2} = 60(c^2-1)^{-5} b_3,\\
		a_{6,3} &= 10(c^2-1)^{-5} b_3,\
		a_{6,4}= (c^2-1)^{-6}b_4\,\
		a_{6,5} = 2(c^2-1)^{-5}b_4,\\
		a_{6,6} &= 45(c^2-1)^{-5} b_3,\
		a_{6,7} = 20(c^2-1)^{-4} b_3,\
		a_{6,8} = (c^2-1)^{-4}b_4,\\
		a_{6,9} &= 15(c^2-1)^{-4} b_3,\
		a_{6,10} = 15(c^2-1)^{-6} b_3,\\
		b_1&= \frac{1}{2160}\left(-3 {\Delta t}^4 c^8-2 {\Delta t}^4 c^6+10 {\Delta t}^2 {\Delta x}^2 c^4-2 {\Delta x}^4 c^2-3 {\Delta x}^4\right),\\
		b_2 &= \frac{1}{720}\left(-2 {\Delta t}^4 c^8-3 {\Delta t}^4 c^6+10 {\Delta t}^2 {\Delta x}^2 c^4-3 {\Delta x}^4 c^2-2 {\Delta x}^4\right),\\
		b_3 &= \frac{1}{302400}\big(10 {\Delta t}^6 c^{12}+22 {\Delta t}^6 c^{10}+3 {\Delta t}^6 c^8-77 {\Delta t}^4 {\Delta x}^2 c^8+28 {\Delta t}^2 {\Delta x}^4 c^6\\
		&\qquad -28 {\Delta t}^4 {\Delta x}^2 c^6-3 {\Delta x}^6 c^4+77 {\Delta t}^2 {\Delta x}^4 c^4-22 {\Delta x}^6 c^2-10 {\Delta x}^6\big),\\
		b_4 &= \frac{1}{120960}\big(72 {\Delta t}^6 c^{12}+94 {\Delta t}^6 c^{10}+9 {\Delta t}^6 c^8-413 {\Delta t}^4 {\Delta x}^2 c^8+112 {\Delta t}^2 {\Delta x}^4 c^6\\
		&\qquad -112 {\Delta t}^4 {\Delta x}^2 c^6-9 {\Delta x}^6 c^4+413 {\Delta t}^2 {\Delta x}^4 c^4-94 {\Delta x}^6 c^2-72 {\Delta x}^6\big).
	\end{align*}
	Refer to \cite{multisymplecticSoftware} for source code.

	\section{Conclusion}
	While backward error analysis of discretisation schemes for PDEs appears intractable, we have introduced a method to perform backward error analysis for {\em symmetric} solutions of numerical methods for PDEs. We show that if a PDE arises from a variational principle and is discretised by a method that respects its variational structure, then symmetric solutions will typically be governed by a modified variational principle up to any order in the discretisation parameters. The modified variational principle is defined by a modified Lagrangian, which is a formal power series in the discretisation parameters. The modified Lagrangian has the correct variational order, i.e.\ order 1, provided that analogous results hold true in the continuous setting.
	Special attention needs to be paid to the fact that in the general case the 1st order modified Lagrangian may only exist in modified variables, blurring the distinction between position and velocity variables.
	If further symmetries are present, the modified Lagrangian can be used to derive formal conserved quantities for the equations that govern the symmetric solutions up to any order in the discretisation parameters.
	The theory is illustrated on rotating travelling waves in the nonlinear wave equation discretised by the multisymplectic five-point stencil.
	
	The analysis performed in this article is purely formal. As future work, it would be interesting to explore settings in which optimal truncation results for modified Lagrangians can be proved, similar to the optimal truncation results for traditional backward error analysis in the ODE setting \cite{GeomIntegration}.

	\section*{Acknowledgments}
	We would like to thank Peter Hydon and Chris Budd for discussions on multisymplecticity and variational principles at the Isaac Newton Institute of the University of Cambridge in 2019.
	This research was supported by the Marsden Fund of the Royal Society Te Ap\={a}rangi and the School of Fundamental Sciences, Massey University, Manawat\={u}, New Zealand.
	The authors would like to thank the Isaac Newton Institute for Mathematical Sciences for support and hospitality during the programme {\em Geometry, compatibility and structure preservation in computational differential equations} (2019) when work on this paper was undertaken. RM would like to thank the
	Simons Foundation for their support. This work was supported by: EPSRC grant number EP/R014604/1.
	The authors acknowledge support from the European Union Horizon 2020 research and innovation programmes under the Marie Skodowska-Curie grant agreement No.~691070.










\medskip
Received for publication January 2022; early access June 2022.
\medskip

\end{document}